\documentclass[reqno]{amsart}
\usepackage{amscd}
\usepackage{amssymb}
\usepackage{graphicx}
\newtheorem{thm}{Theorem}[section]
\newtheorem{dfn}[thm]{Definition}

\newtheorem{prop}[thm]{Proposition}
\newtheorem{cor}[thm]{Corollary}
\newtheorem{lma}[thm]{Lemma}
\newtheorem{lem}[thm]{Lemma}
\newtheorem{eg}[thm]{Example}
\newenvironment{ex}{\begin{eg}\rm}{\end{eg}}
\newtheorem{rem}[thm]{Remark}
\newenvironment{rmk}{\begin{rem}\rm}{\end{rem}}
\newtheorem*{cla}{Claim}
\newenvironment{clm}{\begin{cla}\rm}{\end{cla}}

\newenvironment{pf}{\begin{proof}}{\end{proof}}
\numberwithin{equation}{section}

\newcommand{\R}{{\mathbb{R}}}
\newcommand{\C}{{\mathbb{C}}}

\newcommand{\Z}{{\mathbb{Z}}}

\newcommand{\U}{{\bf{U}}}
\newcommand{\f}{{\bf{f}}}

\newcommand{\E}{{\mathcal{E}}}

\newcommand{\conf}{{\mathcal{C}}}

\newcommand{\A}{{\mathcal{A}}}
\newcommand{\J}{{\mathcal{J}}}

\newcommand{\M}{{\mathcal{M}}}
\newcommand{\Lie}{{\mathcal{L}}}

\newcommand{\Ordo}{{\mathcal{O}}}
\newcommand{\cand}{{\mathcal{W}}}
\newcommand{\tcand}{{\mathcal{V}}}

\newcommand{\sblv}{{\mathcal{H}}}

\newcommand{\End}{\operatorname{End}}
\newcommand{\Hom}{\operatorname{Hom}}

\newcommand{\la}{\langle}
\newcommand{\ra}{\rangle}
\newcommand{\pa}{\partial}
\newcommand{\id}{\operatorname{id}}
\newcommand{\pr}{\operatorname{pr}}

\newcommand{\diag}{\operatorname{Diag}}

\begin{document}

\title{Legendrian Contact Homology in $P \times \R$}

\author{Tobias Ekholm}
\address{Department of Mathematics, University of Southern California, Los Angeles, CA  90089-2532}

\author{John Etnyre}
\address{Department of Mathematics, University of Pennsylvania,
        Philadelphia PA 19105-6395}

\author{Michael Sullivan}
\address{Department of Mathematics, University of Massachusetts,
 Amherst, MA 01003-9305}

%\subjclass{Primary: ; Secondary: }
%\keywords{}

\begin{abstract}
A rigorous foundation for the contact homology of Legendrian submanifolds in a contact manifold of the form $P\times \R$ where $P$ is an exact symplectic manifold is established. The class of such contact manifolds include 1-jet spaces of smooth manifolds. As an application, contact homology is used to provide (smooth) isotopy invariants of submanifolds of $\R^n$ and, more generally, invariants of self transverse immersions into $\R^n$ up to restricted regular homotopies. When $n=3$, this application is the first step in extending and providing a contact geometric underpinning for the new knot invariants of Ng \cite{ng1, ng2, ng3}.
\end{abstract}

\maketitle

%\tableofcontents

\section{Introduction and main results}
Let $P$ be an exact symplectic manifold of dimension $2n$ and let $\theta$ be a primitive of the symplectic form $\omega$ on $P$ (i.e., $\omega=d\theta$). Note that the $1$-form $\alpha=dz-\theta$ on the product manifold $P\times\R$, where $z$ is a coordinate on $\R$, is a contact form (i.e., $\alpha\wedge(d\alpha)^n\ne 0$). The contact structure associated to $\alpha$ is the completely non-integrable hyperplane field $\xi=\ker(\alpha)$ and the Reeb vector field of $\alpha$ is $\pa_z$. We will study contact manifolds $(P\times\R,\xi)$ of this form. (For technical reasons we also require that $P$ has finite geometry at infinity, see Definition \ref{fgi.dfn}.)

An $n$-dimensional submanifold $L\subset P\times\R$ which is everywhere tangent to $\xi$ is called {\em Legendrian} and a continuous $1$-parameter family of Legendrian submanifolds is a {\em Legendrian isotopy}. Legendrian contact homology  was introduced by Eliashberg, Givental, and Hofer in \cite{EGH} and independently, for Legendrian knots in $\R^2\times\R$, by Chekanov \cite{Ch}. It associates a differential graded algebra (DGA) to a Legendrian submanifold $L$, the homology of which is a Legendrian isotopy invariant called the {\em Legendrian contact homology} of $L$. The differential of the DGA is defined by fixing an almost complex structure $J$ on $P$ compatible with $\omega$ and counting boundary-punctured $J$-holomorphic disks in $P$, with boundary on the projection of $L$, which are asymptotic to projections of Reeb chords of $L$ at punctures. In the case when $\dim(P)=2$ the Riemann mapping theorem can be used to reduce such disk counts to combinatorics but in the higher dimensional case, showing that Legendrian contact homology is well-defined, requires analytical work. In \cite{ees2, ees3}, analytical foundations for Legendrian contact homology in the case $P=\R^{2n}=\C^n$ equipped with its standard complex structure were established. It is the purpose of this paper to generalize that result. More precisely, the main result of the paper is the following.

\begin{thm} \label{main.thm}
The contact homology of Legendrian submanifolds of $(P\times\R, \xi)$ is well defined. In particular the stable tame isomorphism class of the DGA associated to a Legendrian submanifold $L$ is independent of the choice of compatible almost complex structure and is invariant under Legendrian isotopies of $L$.
\end{thm}

We prove this theorem by providing the technical details needed to adapt the proof of the corresponding theorem given in \cite{ees2, ees3}. The main difference is related to setting up the functional analytic spaces needed to study spaces of holomorphic disks in the absence of a flat metric on $P$. Also, in \cite{ees2, ees3} the standard complex structure on $\C^n=\R^{2n}$ was used in the definition of the DGA, and the transversality needed to compute the DGA was achieved by perturbation of the Legendrian submanifold. In this paper we show that one can either fix an almost complex structure on $P$ and perturb the Legendrian submanifold or fix the Legendrian submanifold and perturb the almost complex structure in order to achieve the transversality properties needed.

The prime application of Theorem \ref{main.thm} is to the construction of invariants of submanifolds of $\R^n$ via contact geometry. In Section~\ref{knot.sec}, following Arnold \cite{Arnold}, we associate to every submanifold, or more generally, every immersion with transverse multiple points, $M$ of $\R^n,$ its conormal lift, $L_M,$ which is a Legendrian submanifold in the 1-jet space $J^1(S^{n-1})$ of the $(n-1)$-sphere. A smooth isotopy of submanifolds or a regular homotopy through immersions with transverse double points lifts to a Legendrian isotopy. Since $J^1(S^{n-1})$ is an example of one of the contact manifolds to which Theorem \ref{main.thm} applies, we get the following.
\begin{cor}
If $M_1$ and $M_2$ are two isotopic submanifolds of $\R^n,$ then the contact homologies of $L_{M_1}$ and $L_{M_2}$ are isomorphic.
\end{cor}

This construction is the inspiration behind Ng's exciting new invariants of knots in $\R^3,$ \cite{ng1, ng2, ng3}. This paper provides a rigorous contact geometric foundation to the contact homology DGA invariant of knots in $\R^3$ and a subsequent paper will demonstrate that Ng's DGA invariant is indeed the contact homology DGA, and then use this interpretation to extend his work to other dimensions.

As another application, we show that local Legendrian ``knotting'' is global if it is detected by the contact homology DGA. More precisely, if two Legendrian submanifolds of $P\times\R$ are contained in small balls and if their contact homology DGA's (computed locally) differ, then they are not Legendrian isotopic globally in $P\times\R.$ Thus the families of examples of non-Legendrian isotopic Legendrian submanifold in $\R^{2n+1}$ (with its standard contact structure) constructed in \cite{ees1} in combination with Theorem \ref{main.thm} give the following result.
\begin{cor}
In $P \times \R$, there exist infinite families of Legendrian spheres, Legendrian tori, and (when
$P$ is four-dimensional) closed orientable Legendrian surfaces of arbitrary genus which have the
same classical invariants but are pairwise non-Legendrian isotopic.
\end{cor}

\begin{rmk}
Let $P$ be a symplectic manifold with symplectic form $\omega$. It is well known  that the space of almost complex structures on $P$ which are compatible with $\omega$ is contractible. A choice of such an almost complex structure allows us to define Chern classes $c_j(P)$ of the complex bundle $TP$. Most of the results below are stated for the case $c_1(P)=0.$ This is not an essential hypothesis, but it significantly simplifies the expositions and it is satisfied in many applications. Whenever we make use of the $c_1(P)=0$ hypothesis we discuss how this assumption may be removed.
\end{rmk}

The paper is organized as follows. In Section~\ref{def.sec} we give the definition of contact homology for Legendrian submanifolds in $P\times\R$ and prove that it is well-defined by modifying arguments in \cite{ees2, ees3}. The technical details of the necessary modifications are presented in Section~\ref{ban.sec}, where we describe the functional analytic bundles which allow us to view the moduli spaces of holomorphic disks used in the definition of the differential in the DGA as the $0$-set of a section, and in Section \ref{mainann.sec}, where the important features (Fredholm-property, transversality, and orientability of the $0$-set) of this section are studied. Finally, in Section~\ref{app.sec} we discuss applications.

\section{Contact homology}\label{def.sec}
We begin this section by describing the contact manifolds for which we will define contact homology in Subsection~\ref{thecontactmanifold}. The DGA of a Legendrian submanifold is defined as a graded algebra in Subsection~\ref{theDGA}. Its differential is defined, using moduli spaces of holomorphic disks in Subsection~\ref{thedifferential} where we also outline proofs of the main properties of the moduli spaces, postponing details of the necessary constructions to later sections. Finally, in Subsection~\ref{invariance}, we demonstrate that contact homology, and the stable tame isomorphism class of the DGA, is invariant under Legendrian isotopies and that it is independent of the choice of almost complex structure.

\subsection{The contact manifold}\label{thecontactmanifold}
Let $(P,\omega)$ be an exact symplectic $2n$-manifold. That is, $\omega$ is a  symplectic form on $P$ such that $\omega=d\theta$,
where $\theta$ is a $1$-form on $P$. Fix such a form $\theta$ and consider the product manifold $P\times\R$. Define the contact form
$$
\alpha=dz -\theta
$$
on $P \times \R$, where $z$ is a linear coordinate along the $\R$-direction. Note that the Reeb field of the contact form
$\alpha$ is simply ${\pa_z}.$ Let
$$
\Pi_P: P \times \R \rightarrow P
$$
denote the projection which forgets the $\R$-factor.

Let $\J(\omega)\to P$ denote the bundle with fiber over $p\in P$ equal to the space of complex structures $J\colon T_p P\to T_p P$ which are compatible with $\omega$. That is, $J^2=-\id$, $\omega$ is positive on $J$-complex lines, and $J$ is an $\omega$-isomorphism. It is well known that the fibers of this fibration are contractible.

Exact symplectic manifolds cannot be closed. In order to achieve compactness for certain spaces of $J$-holomorphic curves for non-compact $P$, we adapt a notion of Gromov (\cite{g} and Chapter 5 of \cite{a-l}). Let $J$ be an almost complex structure compatible with $\omega$ and let $g_J$ be the Riemannian metric $\omega(\cdot, J \cdot).$ Let $B(p,r) \subset P$ be a ball of radius $r > 0$ centered at $p \in P$ with respect to $g_J.$

\begin{dfn} \label{fgi.dfn}
{\rm{ We say $(P, \omega, J)$ has {\em{finite geometry at infinity}} if $g_J$ is complete and if there exists constants $\rho, C >0$ such that the following hold:
\begin{itemize}
\item for all $p \in P$, the map $\exp_p:B(0,\rho) \rightarrow B(p, \rho)$ is a diffeomorphism;
\item for all $r \le r_0, p \in P$, every loop $\gamma \subset B(p,r)$ bounds a disk in $B(p,r)$ of
area less than $C \mbox{ length}(\gamma)^2.$
\end{itemize}
}}
\end{dfn}
We sometimes say $(P, \omega)$ has finite geometry at infinity if it admits a compatible $J$ such that $(P,\omega,J)$ does. Henceforth in this paper we only consider contact manifolds $P\times \R$ where $(P,\omega)$ has finite geometry at infinity. Note that if $(P,\omega,J)$ has finite geometry at infinity,
and if $J'$ is obtained from $J$ via a compact perturbation, then $(P,\omega,J')$ does as well. A specific example, which motivates this paper, is the following.

\begin{ex} \label{1jet.ex}
Let $M$ be a smooth manifold and let $J^1(M)$ denote its $1$-jet space. Note that $J^1(M)=T^\ast M
\times \R$. Let $\pr_M\colon T^\ast M\to M$ denote the natural projection. The {\em canonical or
Liouville} $1$-form $\lambda_M$ on $T^\ast M$ is
$$
\lambda_M(V)=\beta(d\pr_{M}V)\quad\text{ for } V\in T_\beta(T^\ast M).
$$
The {\em standard symplectic form} on $T^\ast M$ is the $2$-form $\omega_M=-d\lambda_M$ and the
{\em standard contact form} on $J^1(M)$ is the $1$-form $\alpha_M=dz-\lambda_M$. If
$q=(q_1,\dots,q_n)$ are local coordinates on $M$ then $(q,p)=(q_1,p_1,\dots,q_n,p_n)$ are local
coordinates on $T^\ast M$, where $(q,p)$ corresponds to the covector
$$
p_1\,dq_1+\dots+p_n\,dq_n\in T^\ast_q M.
$$
In these local coordinates we have
$$
\lambda_M=p\,dq=\sum_j p_j\,dq_j \quad \text{and}\quad  \omega_M=dq\wedge dp=\sum_j dq_j\wedge dp_j.
$$
It is easy to see that $(T^*M, \omega_M)$ has finite geometry at infinity.
\end{ex}

%----------------------------------------------------------------
\subsection{The graded algebra}\label{theDGA}
%----------------------------------------------------------------
Contact homology associates a graded algebra to a Legendrian submanifold of a contact manifold. The
construction of this algebra for a Legendrian submanifold $L$ in $P \times \R$ is very similar to
the one presented in \cite{ees2} (where $P \times \R$ is standard contact $\R^{2n+1}$). We briefly review it here.

Let $c$ denote a {\em{Reeb chord}} of $L$, that is, a trajectory of the Reeb vector field $\pa_z$
starting and ending at points on $L.$ Then $c^* := \Pi_P(c)$ is a double of $\Pi_{P}(L).$ We say
that $L$ is {\em chord generic} if the only self intersections of the Lagrangian immersion $\Pi_P(L)$ are transverse double points. (Note that this is an open and dense condition.) For chord generic $L,$ let $\A$ be the free associative,
non-commutative algebra over $\Z[H_1(L)]$ generated by the (finite set of) double points in
$\Pi_{P}(L).$ (In fact, we will need slightly stronger genericity conditions on $L$. In the language \cite{ees2} we assume $L$ is generic among admissible
chord generic Legendrian submanifolds.)

Assume $L$ is connected and that $H_1(P)$ is free. In this case the algebra $\A$ is graded with grading in $\Z / c(P,\omega)\Z$. Here, $c(P,\omega)$ is the generator of the image of $2c_1(P, \omega)\colon H_2(P;\Z)\to\Z,$ where $c_1(P,\omega)$ is the first Chern class of the tangent bundle of $P$ equipped with any almost complex structure compatible with $\omega.$ Note that the grading is in $\Z$ when $c(P,\omega)=0,$ as is the case when $P$ is the cotangent bundle of a manifold. When $P$ has torsion in its first homology then the grading on $\A$
can be taken to be a $\Z/2\Z$ grading or a rational grading. We will not discuss this situation but
for details see for example \cite{EGH}. To define the grading on $\A$ we need to fix embedded
circles $h_1,\ldots, h_k$ whose homology classes generate the first homology of $P.$ In addition we
must fix a symplectic trivialization of the tangent bundle $TP$ of $P$ over each $h_i.$ For a Reeb
chord $c$ write $c \cap L = \{ c^-, c^+\}$ where $c^+$ has the larger $z$ coordinate. To define the
grading on a double point $c^*=\Pi_P(c)$ choose a ``capping path'' $\gamma_c$ in $L$ that runs from
$c^+$ to $c^-.$ There is a surface $\Sigma_c$ such that $\partial \Sigma_c
=\Pi_P(\gamma_c)-\sum_{i=1}^k n_i h_i,$ for some integers $n_i.$ (Note the $n_i$ are uniquely
defined given that $H_1(P)$ has no torsion.) There is a unique trivialization of $TP$ over
$\Sigma_c$ that extends the chosen trivializations over $h_i.$ With respect to this trivialization
of the symplectic bundle $TP$ over $\Pi_P(\gamma_c)$ we can think of $\Pi(T\gamma_c L)\subset TP$ as a path of
Lagrangian subspaces in $\C^n.$ Closing this path in the standard way (see \cite{ees2}) we get a
loop in the Grassmanian of Lagrangian subspaces in $\C^n.$ We define the Conley-Zehnder index $\nu(c)$ of $c$
to be the Maslov index of this loop. See Section 1 of \cite{ees2} to review the Maslov and
Conley-Zehnder indices. The grading on $c$ is $|c|=\nu(c)-1$ and the grading of a homology class $A
\in H_1(L)$ is defined as the negative Maslov index of the loop of Lagrangian subspaces in $P$
tangent to $\Pi_P(L)$ along some loop $\beta \subset L$ representing $A.$ (Once again we must choose
a surface realizing the homology between $A$ and our generators $h_i$ to get a loop of Lagrangian subspaces.)
Note the grading on elements of $\A$ depend on the various choices of surfaces $\Sigma_c,$ but two
different choices will define a two-dimensional homology class in $P$ and the gradings will differ
by the value of $2c_1(P,\omega)$ on this homology class. Thus our gradings are defined modulo
$c(P,\omega).$ When $L$ is disconnected, $\A$ is graded over $\Z_2$, see \cite{ees3}.

%----------------------------------------------------------------
\subsection{The differential}\label{thedifferential}
%----------------------------------------------------------------
To define the differential $\pa\colon \A \rightarrow \A$, we need to assume that
$L$ is spin and choose a spin structure $\mathfrak{s}.$ (Without the spin condition we could still
define $\A$ and $\pa$ over $\Z_2.$) We note that the condition that $L$ be spin can be somewhat
weakened \cite{FOOO}, but we will only describe the spin case for the sake of simplicity and the fact
that our applications all satisfy this condition, see Subsection \ref{spin.sec}. This differential will be determined by counting
certain $J$-holomorphic disks.

For convenience we fix some complete Riemannian metric $h$ on $P$ which we will refer to as the {\em
background} metric. When speaking about distances in $P$ we will, unless otherwise explicitly
stated, refer to this metric. Sometimes we speak of distances also in $P \times \R$ in which case we
refer to the metric $h+dz^2$ on $P \times \R$.

An almost complex structure $J$ on $P$, compatible with $\omega$, will be called {\em adapted} to
$L$ if:
\begin{itemize}
\item
For each Reeb chord $c$ of $L$, there is some ball $B(c^\ast, r_c)\subset P$ (measured with the
background metric) around the double point $c^\ast$ with coordinates $x+iy\in\C^n$ so that in these
coordinates $J$ agrees with the standard complex structure in $\C^n$.
\item $(P,\omega,J)$ has finite geometry at infinity.
\end{itemize}

Given such an $L$ and $J$, we say that $L$ is {\em admissible} if there are neighborhoods
$W^\pm\subset L$ around each Reeb chord end point on $L$ such that $\Pi_P(W^\pm)\subset
B(c^\ast,r_c)$ are real analytic submanifolds in the coordinates mentioned above.

Let $\{c_j\}_{j=1}^r$ be the set of Reeb chords of $L$. Let $D_{m+1}$ denote the unit disk in the complex plane with $m+1$ punctures $p_0,\dots,p_m$ on the boundary listed in counter-clockwise order. A {\em $J$-holomorphic disk with boundary on $L$} is a pair of smooth maps $(u,f),$ $u\colon D_m\to P$
and $f\colon \pa D_m\to\R,$ such that
\begin{itemize}
\item $u$ solves the $\bar\pa_J$-equation,
$$
\bar\pa_J u= du + J\circ du\circ i=0.
$$
\item $u$ is asymptotic to Reeb chords at its punctures,
$$
\lim_{\zeta\to p_j}u(\zeta)=c_k^\ast,
$$
for some $k$ and every $j$.
\item The map $(u|\pa D_m, f)\colon \pa D_m\to P \times \R$ has image in $L$.
\end{itemize}

We distinguish positive and negative punctures of a $J$-holomorphic disk. Let $(u,f)$ represent such a disk and let $p\in \pa D_m$ be a puncture. Note that the complex orientation on $D_m$ induces an orientation on $\pa D_m$. If $\Omega\subset \pa D_m$ is an arc and if $p\in\Omega$ then an orienting tangent vector of $\pa D_m$ at $p$ points into one of the components of $\Omega\setminus\{p\}$. We say that points in this component lies in the {\em positive direction} of $p$ and that points in the other component lies in the {\em negative direction} of $p$. We say that a puncture $p$ is positive if $(u,f)$ takes points in $\pa D_m$ in the positive (negative) direction of $p$ to the upper (lower) sheet of $L$. We say that $p$ is negative if $(u,f)$ takes points in $\pa D_m$ in the positive (negative) direction of $p$ to the lower (upper) sheet of $L$.

To simplify terminology, we will often say ``$u$ is a holomorphic disk'' instead of ``$(u,f)$ is a $J$-holomorphic disk with boundary $L.$''

If $a$ is a Reeb chord of $L$ and if ${\mathbf b}:=b_1\ldots b_m$ is a (possibly empty) ordered collection (a word) of Reeb chords then let $\M_A(a;{\mathbf b}) = \M_A(a;{\mathbf b};J,L)$ be the
moduli space of holomorphic disks (up to conformal reparameterization) with the following boundary data. Maps $u$ representing elements in $\M_A(a;{\mathbf b})$ have exactly one positive puncture $p_0$ and $u(p_0) = a^\ast,$ at the ordered $m$ negative punctures $p_1,\dots,p_m$, $u(p_i) = b_i^\ast,$ and the image of the boundary $(u(\pa D_{m+1}),f(\pa D_{m+1}))\subset L$ completed with the ``capping paths'' $\gamma_a$ and the $\gamma_{b_i}$'s
(see the previous subsection) represents the homology class $A\in H_1(L).$

\begin{prop}\label{mfld.prop}
Let $L\subset P\times\R$ be a chord generic Legendrian submanifold and assume that $c(P,\omega)=0$.
\begin{enumerate}
\item Fix a positive integer $N.$ For an open dense set of  $J$ in the space of those almost complex
structures adapted to $L$ for which $L$ is admissible with respect to $J,$ $\M_A(a;{\mathbf b})$ is a manifold provided
\[
|a| - \sum_{i=1}^m |b_i| + \mu(A) -1 \le N.
\]
\item Fix $J$ adapted to $L$ for which $L$ is admissible. Then for any generic perturbation of $L,$
in the space of admissible Legendrian submanifolds to which $J$ remains adapted, $\M_A(a;{\mathbf b})$ is a manifold provided
\[
|a| - \sum_{i=1}^m |b_i| + \mu(A) -1\le 1
\]
where $\mu$ is the Maslov index.
\end{enumerate}
In both cases, the dimension of $\M_A(a;{\mathbf b})$ is $|a| - \sum_{i=1}^m |b_i| + \mu(A) - 1$ and
the manifold is compact in the sense of Gromov (see below). In particular, if the dimension is 0
then $\M_A(a;{\mathbf b})$ is a finite set of points. Moreover, if $L$ is a equipped with a spin structure
$\mathfrak{s}$  and if certain capping operators at its Reeb chords have been specified, see \cite{ees3}, then $\M_A(a;{\mathbf b})$ is naturally oriented.
\end{prop}

\begin{rmk}
In the case when $c(P,\omega)\ne 0$, an analogous lemma holds. The main difference in that case is that the formal dimension of a moduli space containing a holomorphic disk depends on the homology class represented by $u$ in $H_2(P,\Pi_P(L))$. However, the expression for the dimension given in Proposition \ref{mfld.prop} still holds modulo $c(P,\omega)$.
\end{rmk}

\begin{proof}
Let $L\subset P\times\R$ be  Legendrian submanifold and let $J$ be an almost complex structure on $P$ adapted to $L.$ 
In Subsection \ref{sec.coords} we construct a Banach manifold of maps $\cand_{2,\epsilon}(a{\mathbf b})$ the elements of which are triples $(u,f,\kappa)$, where $u\colon D_m\to P$ and $f\colon \pa D_m\to\R$ are such that:
$(u,f)|\pa D_m$ maps into $L;$ 
$u(p_0) = a^*$ and $u(p_i) = b_i^*;$
and the restriction of $\bar\pa_{J,j_\kappa} u$ to $\pa D_m$ 
(see \eqref{barpasec.eq}) equals $0.$  
Here $\kappa$ is a conformal structure on $D_m$ and $j_\kappa$ a corresponding complex structure. In Subsection \ref{sec:bundle} we construct a bundle $\E=\E_{1,\epsilon}$ over $\cand_{2,\epsilon}(a{\mathbf b})$ the fiber of which over $(u,f,\kappa)$ consists of $(J,j_\kappa)$-complex anti-linear maps $T D_m\to u^\ast(TP)$. Thus, the $\bar\pa_{J,j_\kappa}$-operator
\begin{equation}\label{barpasec.eq}
\bar\pa_{J,j_\kappa}u=du+J\circ du\circ j_\kappa
\end{equation}
gives a section $\cand_{2,\epsilon}\to \E$ and the $J$-holomorphic disks correspond to zeroes of this section.

Let $\Lambda$ be a space of compatible almost complex structures (for the proof of (1)) or a space of Legendrian submanifolds (for the proof of (2)). We show in Subsections  \ref{DLs.sec} and \ref{DJ.sec}, respectively how to patch the bundles discussed above to a bundle $\E_\Lambda\to\cand_{2,\epsilon,\Lambda}$, where $\cand_{2,\epsilon,\Lambda}\to\Lambda$ is a bundle the fiber of which over $\lambda\in\Lambda$ is $\cand_{2,\epsilon}$ as defined using the almost complex structure or Legendrian submanifold corresponding to $\lambda$. In Subsection \ref{linearizeJ.sec} we compute the linearization of the section $\bar\pa\colon\cand_{2,\epsilon,\Lambda}\to\E_\lambda$, where $\bar\pa$ is given by the expression \eqref{barpasec.eq} at a point $((u,f,\kappa),\lambda)$.  In Subsection \ref{fred.sec} we demonstrate that the linearization $L\bar\pa_\lambda$, for fixed $\lambda$, of this section is a Fredholm operator and we compute its index in Lemma \ref{fred.lma}. Finally, in Subsection \ref{trans.sec}, we prove that the full linearization (linearizing also with respect to the $\lambda$-variable) of the section $\bar\pa$ is surjective if $\Lambda$ contains all admissible variations of the complex structure  or  the Legendrian submanifold, see Lemma \ref{trans.lma}. It is at this point that the proofs for (1) and (2) differ. To prove the second statement we establish transversality on the complement of exceptional holomorphic disks in $\cand_{2,\epsilon,\Lambda}$. (Exceptional holomorphic disks are defined in Definition 6.10  of \cite{ees2}.) With these transversality properties established we prove (1) and (2) for a Baire set of almost complex structures and Legendrian submanifolds, respectively,  as follows.

Let ${\mathcal Z} = \bar\pa^{-1}(0)$. Since the linearization of $\bar\pa$ is surjective ${\mathcal Z}\subset \cand_{2,\epsilon,\Lambda}$ is a submanifold. Moreover, the map $\pi\colon{\mathcal Z}\to \Lambda$ is Fredholm with index equal to the index of $L\bar\pa_\lambda$. An application of the Sard-Smale theorem then shows that for generic $\lambda$, $\pi^{-1}(\lambda)\cap{\mathcal Z}$ is a submanifold of dimension equal to the Fredholm index of $L\bar\pa_\lambda$. This shows that (1) holds. In case (2) an application of the argument given in the proof of Theorem 6.15 of \cite{ees2} rules out the existence of exceptional disks of low dimensions and completes the proof of (2). (Since the proof of the first statement does not require special considerations for such disks, see the proof of Lemma \ref{trans.lma}, the dimension can be arbitrarily high in this case.)

We will not prove or explain the Gromov compactness property in this paper since the discussion is identical to the $P \times \R = \R^{2n+1}$ case given in Sections 1 and 8 of \cite{ees2}. (Note that Lemma 8.3 of \cite{ees2} carries through because locally at any double point of $\Pi_P(L)$, $(P,J)$ looks like $\C^n.$) For the reader familiar with Gromov compactness in other situations, we point out three relevant facts:
\begin{itemize}
\item Because $\Pi_P(L)$ is an exact Lagrangian immersion in an exact symplectic manifold, there are
neither non-constant smooth (at the boundary) holomorphic disks with boundary on $\Pi_P(L)$ nor non-constant
holomorphic maps of closed Riemann surfaces to $P$.
\item Any disk in $\M_A(a;{\mathbf b})$ has ($\omega$-)area equal to $z(a^+) - z(a^-) - \sum_i
(z(b^+_i) - z(b^-_i))$ where $z$ refers to the last coordinate in $P \times \R.$
\item The above area property together with finite geometry at infinity imply that for any set of adapted almost complex structures $J_\lambda,$ $\lambda \in \Lambda$, which agree outside some compact subset of $P$, there exists a compact set $K\subset P$ independent of $\lambda$ such that all $J_\lambda$-holomorphic disks with boundary on $L$ lie inside $K.$
\end{itemize}

Gromov compactness implies that the moduli spaces are compact. Since the transversality conditions used above are open it follows that (1) and (2) hold for open dense subsets of the space of almost complex structures and Legendrian submanifolds, respectively.

Finally, the orientation of the manifolds are induced by equipping the determinant bundles of the linearization of the $\bar\pa$-operators with orientations in a specific manner. The details of this construction are discussed in Subsection \ref{spin.sec}.
\end{proof}

We can now define the differential for a generic admissible Legendrian submanifold. Let $\mathcal{C}$ be its set of Reeb chords, and let $\mathcal{A}$ denote its algebra. For any generator $a\in \mathcal{C}$ of $\mathcal{A}$ we set
\begin{equation}
        \partial a =\sum_{\hbox{dim } \mathcal{M}_A(a;{\mathbf b})=0} (-1)^{(n-1)(|a|+1)}
        \bigl(\#\mathcal{M}_A(a;{\mathbf b})\bigr)A {\mathbf b},
\end{equation}
where $\# \mathcal{M}$ is the number of points in $\mathcal{M}$ counted with signs induced by the
orientation, and where the sum ranges over all words ${\mathbf b}$ in the alphabet $\mathcal{C}$ and
$A\in H_1(L)$ for which the above moduli space has dimension 0. Note that by the second itemized
point in the above proof, we can conclude that the sum is finite and $a$ is not a letter in
${\mathbf b}.$ We then extend $\partial$ to a map $\partial:\mathcal{A}\to \mathcal{A}$ by linearity
and the graded Leibniz rule.

\begin{lem}\label{dsquared=0.lem}
        $\partial\circ \partial=0$
\end{lem}

\begin{proof}
This is a standard argument following from Gromov compactness and a gluing theorem, Proposition~\ref{glud^2=0.prop}. It is also necessary to check that the orientations have the correct behavior under gluing, see Section \ref{spin.sec}.
\end{proof}

We refer to the pair $(\A, \pa)$ as the {\em{DGA}}  associated to $L$
and $J.$ The {\em{contact homology}} of $L$ and $J$ is
\begin{equation} \notag
        HC_*(L,J)= \hbox{Ker } \partial/ \hbox{Im } \partial.
\end{equation}

%----------------------------------------------------------------
\subsection{Invariance of contact homology}\label{invariance}
%----------------------------------------------------------------
\begin{prop}\label{invariance.prop}
Let $L_s\subset P\times\R$, $s\in[0,1]$ be a Legendrian isotopy such that $L_0$ and $L_1$ are
admissible with respect to adapted almost complex structures $J_0$ and $J_1$, respectively. Then the
DGAs of $(L_0,J_0)$ and $(L_1,J_1)$ are stable tame isomorphic. In particular, $HC_*(L_0,J_0)$ and
$HC_*(L_1,J_1)$ are isomorphic as graded algebras.
\end{prop}

See \cite{ees1} for a review of stable tame isomorphisms. We prove Proposition \ref{invariance.prop} below. To conclude from it that contact homology is a well-defined invariant for Legendrian isotopy classes, we note the following.

\begin{lma}\label{lmaadm0}
Let $L\subset P\times\R$ be any Legendrian submanifold. Then there exists an arbitrarily small
Legendrian isotopy $L_t$ of $L=L_0$ such that $L_1$ is chord generic and admissible with respect to
some almost complex structure $J$ on $P$ adapted to $L$.
\end{lma}

\begin{pf}
For a fixed chord generic Legendrian submanifold the existence of an adapted complex structure is easy to establish using the contractibility of the space of almost complex structures. As in the proof of Lemma 2.5 in \cite{ees2} one can make a Legendrian submanifold real analytic close to its double points.
\end{pf}

Our proof of Proposition \ref{invariance.prop} involves studying how bifurcation moments in a Legendrian isotopy affect the DGA and we will use the cobordism method of \cite{ees3} to study these moments. In order to isolate the bifurcation moments, we define chord genericity and admissibility for $1$-parameter families of Legendrian submanifolds. Essentially, we require the isotopy to be chord generic except for isolated ``birth-death'' moments of pairs of chords (quadratic self-tangencies). We also require that our changing almost complex structure remain integrable in uniform-size neighborhoods, with respect to the background metric, of the double points. Since these are all local considerations the arguments in Sections 1.5 and 2 of \cite{ees2} suffice in the current situation as well. See \cite{ees2} for details.

For technical reasons, see Remark \ref{commute.rmk}, we restrict the isotopies further so that at any point in time either the almost complex structure or the Legendrian is changing, but not both. The following lemma shows that we can do so without loss of generality.

\begin{lma}\label{lmaadm1}
Let $L_s\subset P\times\R$, $s\in[0,1]$ be a Legendrian isotopy such that $L_0$ and $L_1$ are
admissible with respect to adapted almost complex structures $J_0$ and $J_1$, respectively. Then
there exists a chord generic Legendrian isotopy $L'_s$, $s\in[0,1]$, arbitrarily $C^0$-close to
$L_s$, with $L'_0=L_0$ and $L'_1=L_1$ and a family of almost complex structures $J_s$, $s\in[0,1]$
with the following property:
there exists a partition $0=s_0<s_1<\dots<s_N=1$ of $[0,1]$ such that $L_s$, $s\in[s_j,s_{j+1}]$ is
adapted to $J_{s_j}$ and such that $L_{s_{j+1}}$ is adapted to $J_s$ for $s\in [s_j,s_{j+1}]$.
\end{lma}

\begin{pf}
Assume first that the isotopy is chord generic and does not have any self tangencies. In this
case the families of Reeb chords $c^\ast_s$ of $L_s$ are continuous. By compactness of $[0,1]$ there
exists $r_0>0$ such that for every $s\in[0,1]$, every Reeb chord $c_s$, and every $r<r_0$,
$\Pi_P(L_s)\cap B_{r_0}(c^\ast_s)$ are two transversely meeting Lagrangian disks. Let $r<\frac12
r_0$. Again by compactness of $[0,1]$ we find a finite partition $0=s_0<s_1<\dots<s_N=1$ such that
$c_s^\ast\in B_{\frac12 r}(c_{s_j}^\ast)$ for $s\in[s_{j-1},s_{j+1}]$.

We define a path of complex structures inductively. Assume that $J_s$ has been defined for
$s\in[0,s_j]$ in such a way that %$J_{s_k}|B_{r}(c^\ast_{s_k})$ and
$J_{s}|B_{\frac18r}(c^\ast_{s_{k+1}})=J_{s_k}$ for $s\in[s_k,s_{k+1}]$ and $k<j$. We then find
$J_s$, $s_j\le s\le s_{j+1}$. We take $J_{s_{j+1}}$ to be the integrable complex structure over
$B_{r}(c_{s_{j+1}}^\ast)$ which is the conjugation of $J_{s_k}$ by the differential of the scaling
$x\mapsto \frac18 x$. Since the fiber of the bundle of almost complex structures compatible with
$\omega_M$ is contractible we can join $J_{s_k}$ to a complex structure $J_{s_{k+1}}$ agreeing with
the one defined on $B_{r}(c^\ast_{s_{k+1}})$ and so that $J_s$ is fixed on $B_{\frac18
r}(c^\ast_{s_{k+1}})$. By induction we can continue to $j=N$. With this accomplished we use the
argument from \cite{ees2} (approximation of smooth functions with a finite part of its Taylor series) to
make $L_s$ admissible as desired.

In presence of self tangencies we deal with neighborhoods of the self tangency points separately and
apply the above argument to the rest of the isotopy.
\end{pf}

\begin{proof}[Proof of Proposition \ref{invariance.prop}]
Like in \cite{ees2,ees3}, there are two bifurcation moments to consider:
the birth-death of two Reeb chords mentioned above and the isolated appearance
of a non-generic ``handle-slide'' disks. The behavior of holomorphic disks at such a degenerate moment is governed by a result similar to Proposition \ref{mfld.prop} for $1$-parameter families of chord generic Legendrian submanifolds and for Legendrian submanifolds with exactly one degenerate Reeb chord corresponding to a self-tangency double points. In the case $P=\R^{2n}$ such results were proved in \cite{ees2} by a modification of the argument giving the counterpart of Proposition \ref{mfld.prop}. The results necessary in the more general setting can be obtained from the discussion in Section \ref{sec.coords} by a completely analogous modification.

We adapt the cobordism method used in \cite{ees3}. Assume the bifurcation moment occurs at $t=0$ and  for $-1 \le t \le 1 $  no other bifurcation occurs. Let $L_t = \Phi_t(q)=(p(q,t),z(q,t)), q\in L$ be the isotopy and $J_t$ be the family of almost complex structures.

We construct a Legendrian embedding of $\R\times L$ into $(P\times\R^2)\times\R$ with contact form $dz - \theta - ydx$ as follows
\begin{equation}
\notag
(t, q)\mapsto (t, y(q,t), p(q,t), z(q,t))
\quad \mbox{where} \quad
y(q,t)=-\theta\left(\frac{\pa p}{\pa t}\right)+\frac{\pa z}{\pa t}.
\end{equation}
Note that $(P \times \R^2, d\theta + dx\wedge dy, J := J_t \oplus i)$ is an exact symplectic manifold with finite geometry at infinity, where $i$ is the standard complex structure on $\R^2$ and $J_{t \le -1} := J_{-1}$, $J_{t \ge 1} := J_1.$

Note that by Lemma \ref{lmaadm1}, either $J_t$ or $\Phi_t$ is independent of $t.$
Consider first the case when $J_t$ is independent of $t.$ In this case, we repeat the argument in \cite{ees3} which deduces the stable tame isomorphism invariance of the DGA's of $L_0$ and $L_1$ from the fact that the contact homology differential $\Delta$ of $L\times\R$ satisfies $\Delta\circ\Delta=0.$ In the case that $J_t$ varies with $t$ and there is a handle slide disk for $J_0$ the exact same algebraic argument as in the case when $\Phi_t$ varies and there is a handle slide disk for $\Phi_0$ applies, after the counterparts of Lemmas 4.18, 4.19, and 4.20 have been established. This is straightforward, as we illustrate by giving the modification of Lemma 4.18 below.
\end{proof}

Consider the case when a handle-slide disk appears at $J_0$ as $J_t$, $-1\le t\le t$ varies. Like in \cite{ees3}, we perturb the embedded $\R \times L$ by a sequence of Morse functions $f_k$ which converge to the constant function $1$,
each with 2 local minima at $\pm1$ and a local maximum at $0.$
We relabel the perturbed Legendrian $L^k.$
For each Reeb chord $c$ in $L,$ there are three in  $L^k$, $c[\pm1], c[0]$,
with $|c[0]| = |c[\pm1]| +1  = |c| + 1.$
Given $ 0 \le \delta <1$, we also rescale a shrinking interval of $J_t$, $|t| \le \delta$,
to create a new $1$-parameter family of almost complex
structures $J_t^\delta,$ where  $J_{\ge 1}^\delta =J_{\delta}$,
 $J_{\le -1}^\delta =J_{-\delta}$, $J_0^\delta = J_0.$
 Let $J^\delta = J_t^\delta \oplus i.$

The analog of Lemma 4.18 \cite{ees3}, the first step in the proof handle-slide invariance, is the following. (The remaining steps have similar analogues).

\begin{lma}
There exists $k_0$ such that for all $k>k_0$ there exists a $\delta_k>0$
such that for all $\delta<\delta_k$ and any Reeb chord $c$ of $L$ the
following holds. The moduli spaces $\M(c[0],c[1])$ and
$\M(c[0],c[-1])$ of $J^\delta$-holomorphic
disks with boundary on $L^k$ consists of exactly one point
which is a transversely cut out rigid disk.
Moreover the sign of the rigid disk in $\M(c[0],c[1])$ and that of the
disk in $\M(c[0],c[-1])$ are opposite.
\end{lma}

\begin{proof}
When $\delta = 0$, we have $(J^\delta, J_0)$- and $(J^\delta, i)$-holomorphic projections in the two components of $P \times \R^2.$ So the analysis of the moduli spaces is identical to that done in the $\delta=0$ case of the proof of Lemma 4.18 \cite{ees3}. Gromov compactness then implies the result holds for small $\delta.$
\end{proof}

%********************************************************
\section{Functional analytic setup} \label{ban.sec}
%********************************************************
This section provides the necessary functional analytic setup to describe the moduli spaces used in the definition of the contact homology DGA. We describe a Banach manifold of candidate maps ("the configuration space") in Subsection~\ref{sec.coords}. In Subsection~\ref{sec:bundle}, we describe a bundle over this Banach manifold where a suitable $\bar\pa$-operator takes values. The $0$-set of the $\bar\pa$-section will describe the moduli spaces of holomorphic disks. We discuss how to extend this bundle over spaces parameterizing deformations of the Legendrian submanifolds and of almost complex structures in Subsections~\ref{DLs.sec} and~\ref{DJ.sec}, respectively. Finally, in Subsection~\ref{linearizeJ.sec} we discuss the linearization of the section that defines our moduli spaces.

In short, this section provides the setup for the proof of Proposition \ref{mfld.prop}.

%-------------------------------------------------------------------
\subsection{The manifold of candidate maps}\label{sec.coords}
%-------------------------------------------------------------------
Let $L\subset P\times\R$ be a Legendrian submanifold. Fix a complex structure $J$ on $P$ compatible with $\omega$ and adapted to $L$. We first describe a Banach manifold of maps which contain all holomorphic disks with prescribed punctures. Consider $D_m$ with a conformal structure $\kappa$. We fix as in \cite{ees2} a metric on $D_m$ such that neighborhoods of the punctures look like half infinite strips, with coordinates $\tau+it\in\R_\pm+i[0,1]$. We define for $\epsilon\in\R$ weighted Sobolev spaces $\sblv_{k,\epsilon}(D_m,\R^N)$ of functions with $k$ derivatives in $L^2$, by using weights which look like $\exp(\epsilon|\tau|)$ in neighborhoods of the punctures. We denote the complex structure on $D_m$ induced by this metric $j_\kappa$.

\subsubsection{Maps into $P$}
Pick an isometric (with respect to the background metric) embedding $P\subset\R^N$ for some sufficiently large $N$. Fix Reeb chords $b_1,\dots,b_m$ and write ${\mathbf b}=b_1,\dots,b_m$. Fix a reference function $f_{\mathbf b}\colon D_m\to\R^N$ which is a smooth map constantly equal to $b_j^\ast\in\Pi_P(L)\subset P\subset\R^N$ in a neighborhood of $p_j$ for each $j$. Define $\hat\tcand_{2,\epsilon}({\mathbf b};\kappa)$ to be the affine Banach space of functions $u\colon D_m\to\R^N$ such that
\begin{equation}
u-f_{\mathbf b}\in\sblv_{2,\epsilon}(D_m,\R^N).
\end{equation}
Note that elements in $\hat\tcand_{2,\epsilon}({\mathbf b};\kappa)$ are continuous and define $\tcand_{2,\epsilon}({\mathbf b};\kappa)$ as the subset of elements in $\hat\tcand_{2,\epsilon}({\mathbf b};\kappa)$ such that $u(D_m)\subset P$.

\begin{lma}
The subset $\tcand_{2,\epsilon}({\mathbf b};\kappa)$ is closed submanifold of $\hat\tcand_{2,\epsilon}({\mathbf b};\kappa)$. The tangent space $T_{u}\tcand_{2,\epsilon}({\mathbf b};\kappa)$ is canonically isomorphic to the space $\sblv_{2,\epsilon}(D_m,u^\ast(TP))$.
\end{lma}

\begin{pf}
The subspace $\tcand_{2,\epsilon}({\mathbf b};\kappa)$ is closed subset since the $2$-norm controls
the supremum norm. The fact that it is a submanifold and the statement on its tangent space can be
seen as follows. Fix a metric $\eta$ on $\R^N$ in which $P$ is totally geodesic. Let $\exp$
denote the exponential map with respect to this metric. Define the local coordinate map $\Phi\colon\sblv_{2,\epsilon}(D_m,u^\ast(TP))\to \hat\tcand_{2,\epsilon}({\mathbf b};\kappa)$ by letting $\Phi(v)$ be the map
$$
\zeta\mapsto \exp_{u(\zeta)}(v(\zeta)),\quad \zeta\in D_m.
$$
Since the $2$-norm controls the supremum norm and since the injectivity radius of $\eta$ is bounded from below it follows from standard estimates that these maps give local charts with coordinate changes which are $C^1$.
\end{pf}

\subsubsection{Boundary conditions}
Let $L\subset P\times\R$ be a chord generic Legendrian submanifold which is admissible in an adapted
almost complex structure $J$. Consider the Banach manifold
$$
X=\tcand_{2,\epsilon}
({\mathbf  b};\kappa)\times\sblv_{\frac32,\epsilon}(\pa D_m,\R).
$$
Define $\cand_{2,\epsilon}({\mathbf b};\kappa,J)$ to be the closed submanifold of $X$ which consists
of functions $(u,f)$ fulfilling the following conditions
\begin{align}
\label{boundary1.eq}
(u,f)(\zeta)\in L\quad, &\text{ for all }\zeta\in\pa D_m.\\
\label{boundary2.eq}
\int_{\pa D_m}\la\bar\pa_J u,v\ra\, ds=0, &\text{ for all }v\in
C^0_0(\pa D_m, {T^\ast}^{0,1}D_m\otimes u^\ast(TP)).
\end{align}

In order to show that $\cand_{2,\epsilon}({\mathbf b},\kappa,J)$ is a Banach manifold we construct
local coordinates on it using a special $1$-parameter family of metrics on $P$. Consider the metric
$\hat g$ on $T L$ constructed from a metric $g$ on $L$ in \cite{ees2}. Pick an immersion $S\colon T L\to
P$ such that $dS\circ J_L=J\circ dS$ along $L\subset T L$, where $J_L$ is the canonical complex
structure along the $0$-section $L\subset TL$.
Let $U \subset TL$ be a neighborhood of the $0$-section.
Let $B(c_j^\ast,r)\subset P$ be a ball of radius $r$ in the background metric on $P$.
Note that for $r>0$ sufficiently small
$\Pi_P^{-1}(B_r(c_j^\ast))=B_j^+(r)\cup B_j^-(r)$ where $B_j^\pm(r)$ are disjoint (topological) balls in $L$
around the upper and lower endpoints of $c_j$. Fix a small $\delta>0$ and choose a $1$-parameter
family of metrics $g^\sigma$ on $P$ such that
\begin{itemize}
\item $g^\sigma=\hat g$ in $S(U|L\setminus\cup_j(B_j^+(10\delta)\cup
B_j^-(10\delta))$, $0\le \sigma\le 1$.
\item
$g^\sigma=\hat g$ in $S(U|B_j^-(5\delta))$ for $\sigma$ in a
neighborhood of $0$.
\item
$g^\sigma=\hat g$ in $S(U|B_j^+(5\delta))$ for $\sigma$ in a
neighborhood of $1$.
\end{itemize}

Next fix a smooth function $\rho\colon P\times\R\to[0,1]$ such that $\rho|(L\setminus \cup_j
B^+_j(9\delta))=0$ and such that $\rho |B^+_j(6\delta)=1$. Note that if $\gamma\colon [0,1]\to L$ is a
curve then there exists $\eta$ such that the metric $g^{\rho(\gamma(t))}$ on the ball
$B_\eta(\gamma(t))\subset P$ agrees with $\hat g$ for all $t$.

For $(u,f)\in X$, pick an extension $F\in\sblv_2(D_m,\R)$ of $f$.

\begin{lma}\label{lma.coords}
The subspace $\cand_{2,\epsilon}({\mathbf b};\kappa)$ is a smooth submanifold of $X$. Its tangent
space at $(u,f)$ is isomorphic to the closed subspace of all $v\in\sblv_{2,\epsilon}(D_m,u^\ast(T
P))$ which satisfies
\begin{align}
v(\zeta)\in T_{u(\zeta)}L\subset T_{u(\zeta)}P &\text{ for all }\zeta\in\pa D_m,\\ \int_{\pa D_m}\la
\bar\nabla_J v,w\ra\,ds=0 &\text{ for all }w\in C^0_0(\pa D_m, {T^\ast}^{0,1}D_m\otimes u^\ast(TP)),
\end{align}
where
$$
\bar\nabla_J v = \nabla v + J\circ \nabla v\circ j_\kappa,
$$
where $\nabla$ is the connection on $u^\ast (TP)$ induced from the metrics $g^{\rho(F(\zeta))}$.
Moreover, the map
$$
\Psi\colon T_{(u,f)}\cand_{2,\epsilon}({\mathbf b};\kappa)\to\cand_{2,\epsilon}({\mathbf b};\kappa),
$$
$$
\Psi(v)(\zeta)=\exp_{u(\zeta)}^{\rho(F(\zeta))}(v(\zeta)),
$$
gives local $C^1$-coordinates on an $\eta$-ball around $(u,f)$ where $\eta$ depends continuously on
the minimal injectivity radius of the metrics $g^\sigma$, $0\le\sigma\le 1$.
\end{lma}

\begin{pf}
It is straightforward to check that $\cand_{2,\epsilon}({\mathbf b};\kappa)$ is a closed subset. Let
$\exp^\sigma$ denote the exponential map in the metric $g^\sigma$ on $P$. Noting that for
$(u,f)\in\cand_{2,\epsilon}({\mathbf b};\kappa)$, the metric $g^{\rho(f(\zeta))}$ is constantly
equal to $\hat g$ on $T^\ast L$ for $\zeta\in\pa D_m$, standard arguments show that
$$
v\mapsto \exp_{u(\zeta)}^{\rho(F(\zeta))}(v(\zeta))
$$
gives local coordinates around $u$ in $\tcand({\mathbf b};\kappa)$.

We need to check the boundary conditions. The first boundary condition is clearly met if and only if
$v(\zeta)$ is tangent to $L$ for $\zeta$ in the boundary since a neighborhood of $u(\zeta)$ in $L$
is totally geodesic in the metric $g^{\rho(F(\zeta))}$ for all $\zeta\in\pa D_m$.

To see that the second condition is met if and only if $\bar\nabla_J v=0$ on the boundary we compute
$\bar\pa_J(\exp_{u(\zeta)}(v(\zeta)))$ (we drop the superscript on the metric since it is constant
for $\zeta\in\pa D_m$). Note first that in complex coordinates $\tau+it$ on $D_m$ we have
$$
\bar\pa_J w=\frac{\pa w}{\pa \tau}+ J\frac{\pa w}{\pa t}.
$$
Now as in \cite{ees2}
\begin{align*}
&\frac{\pa}{\pa \tau}\exp_{u(\zeta)}(v(\zeta))=X_1(1),\\
&\frac{\pa}{\pa t}\exp_{u(\zeta)}(v(\zeta))=X_2(1),
\end{align*}
where $X_1$ and $X_2$ are the Jacobi fields along the geodesic $s\mapsto \exp_{u(\zeta)}(sv(\zeta))$
with
\begin{align*}
X_1(0)&=\frac{\pa u}{\pa\tau},\\
\nabla_s X_1(0)&=\nabla_\tau v,\\
X_2(0)&=\frac{\pa u}{\pa t},\\
\nabla_s X_2(0)&=\nabla_t v.
\end{align*}
Recall, see Lemma 4.6 in \cite{ees2}, that $X$ is a Jacobi field along a geodesic $\gamma\subset L$
in the metric $\hat g$ if and only if $JX$ is a Jacobi field. It follows that
$\bar\pa_J(\exp_{u(\zeta)}(v(\zeta)))$ vanishes along the boundary if and only if $\bar\nabla_Jv=0$.
\end{pf}

\subsubsection{Conformal structures}\label{confstr.sssec}
The space of conformal structures on a $m$ punctured disk with one distinguished puncture (the positive puncture in our case) is an $(m-3)$-dimensional simplex. We will denote this space $\conf_m$. More concretely, we think of this space as follows. Fix the distinguished puncture at $1\in\pa D$, where $D$ is the unit disk in the complex plane, and the two punctures immediately following it at $i$ and $-1$. Then the positions of the remaining punctures in the lower half of $\pa D$ give coordinates on $\conf_m$. We write
$$
\cand_{2,\epsilon}
({\mathbf b})=
\bigcup_{\kappa \in \conf_m}\cand_{2,\epsilon}(\mathbf b;\kappa).
$$
This is a locally trivial bundle as in \cite{ees2}. (We precompose with diffeomorphisms of $D_m$ which are holomorphic on the boundary.) Moreover, we may view the tangent space $T_\kappa\conf_m$ as spanned by a finite number of elements in $\End^{0,1}(TD_m)$ (complex anti-linear endomorphisms of $TD_m$) of the form $\bar\pa v$, where $v$ is a vector field generating a diffeomorphism which moves one of the punctures.

%----------------------------------------------------------------------
\subsection{Bundles of complex anti-linear maps}\label{sec:bundle}
%----------------------------------------------------------------------
This subsection constructs the bundle of complex anti-linear maps over $\cand_{2,\epsilon} ({\mathbf b}).$ We begin by discussing the bundle of complex anti-linear maps in general and then turn attention to the specific bundles which will be used below.
\subsubsection{Complex anti-linear maps}
Let $B$ be a manifold and let $E$ and $F$ be complex vector bundles over $B$. We think of these
complex vector bundles as real vector bundles with complex structures $J_E$ and $J_F$ on them. Let
$J_E^\lambda$, $\lambda\in\Lambda$ and $J_F^\gamma$, $\gamma\in\Gamma$ be smooth families of complex
structures on $E$ and $F$. A real linear bundle homomorphism $A\colon E\to F$ is
$(J_E^\lambda,J_F^\gamma)$-complex anti-linear if
$$
A - J_F^\gamma \circ A \circ J_E^\lambda=0.
$$
Form the bundle
$$
\Hom^{0,1}(E,F)\to \Lambda\times\Gamma
$$
of $(J^\lambda_E,J^\gamma_F)$-complex anti-linear bundle homomorphisms from $E$ to $F$. We will use
the following local trivialization of this bundle. If $0\in\Lambda$ and $0\in\Gamma$ then each
complex structure in a neighborhood of $J^0_E$ and $J^0_F$ can be written as
\begin{align*}
J^\lambda_E &= J^0_E(1+S^\lambda_E)(1-S^\lambda_E)^{-1},\\
J^\gamma_F &= J^0_F(1+S^\gamma_F)(1-S^\gamma_F)^{-1},
\end{align*}
where $S^\lambda_E=(J_E^\lambda+J^0_E)^{-1}(J_E^\lambda-J^0_E)$ (note that $J_E^0
S^\lambda_E+S^\lambda_EJ_E^0=0$) and similarly for $S^\gamma_F$.

It is straightforward to check that the transformation of bundle homomorphisms
$$
A^0\mapsto (1-S^\gamma_F) A^0 (1+ S^\lambda_E),
$$
takes $(J^0_E,J^0_F)$-complex anti-linear maps to $(J^\lambda_E,J^\gamma_F)$-complex anti-linear maps.
We use such transformations as local trivialization of the (anti)-complex vector bundle $\Hom^{0,1}(E,F)$.
This construction will be used repeatedly below.

\subsubsection{The bundle where the $\bar\pa_J$-operator takes values}
Let $L\subset P\times\R$ be a chord generic admissible Legendrian submanifold with respect to an
adapted almost complex structure $J$ on $P$. We consider the bundle $\E = \E_{1,\epsilon}$ over
$\cand_{2,\epsilon}({\mathbf b})$ the fiber of which over $(u,f,\kappa)$ consists of $(j_\kappa,J)$-complex anti-linear maps
$$
T^\ast D_m\to u^\ast(TP),
$$
which lie in the natural Sobolev space.
In future subsections, we will consider different Legendrian submanifolds, $L_\lambda$ or
almost complex structures, $J_\lambda.$
To indicate this, we will relabel the bundle as $\E(L_\lambda)$
or $\E(J_\lambda).$
We find local trivializations of this bundle as follows: let
$U\times K$ be an open neighborhood of
$(0,\kappa)$ in $\sblv_{2,\epsilon}(u^\ast(TP))\times \conf_m$.
Let $\Phi$ denote the local coordinate map from $U$ to $\cand_{2,\epsilon}({\mathbf{b}})$.
Let $(v, \mu) \in U \times K.$
Note that parallel-translation in the metric $g^{\rho(F(\zeta))}$ along the unique geodesic
from $u(\zeta)$ to $\Phi(v)(\zeta)$ identifies $(J,j_\kappa)$-complex anti-linear maps
$$
T^\ast D_m\to (\Phi(v))^\ast(TP)
$$
with $(J_v, j_\kappa)$-complex anti-linear maps
$$
T^\ast D_m\to u^\ast (TP).
$$
Here $J_v=\Pi^{-1}_v\circ J\circ \Pi_v$, where $\Pi_v$ denotes parallel translation.
We can thus
trivialize the bundle as in the previous subsection:
$$
A\mapsto (1-S_v)A(1+\gamma),
$$
where $S_v=(J+J_v)^{-1}(J_v-J)$ and $\gamma=(j_\kappa+j_\mu)^{-1}(j_\mu-j_\kappa)$. Note also that the
$\bar\pa_{J,\kappa}$-operator gives a smooth section
$$
\bar\pa_{J,\kappa}\colon \cand_{2,\epsilon}({\mathbf b})\to\E.
$$

\subsection{Deforming the Legendrian submanifold}
\label{DLs.sec}

Let $L_\lambda$, $\lambda\in\Lambda$ be a family of chord generic Legendrian submanifolds where
$\Lambda$ is an open ball in some Banach space. Let $0\in\Lambda$ and assume that $J$ is an almost
complex structure adapted to $L_0$ and that $L_0$ is admissible with respect to $J$. Note that by
continuity there exists some neighborhood $U$ of $0$ in $\Lambda$ so that $L_\lambda$ is adapted to
$J$ for all $\lambda\in U$. We will restrict attention to deformations such that for all $\lambda\in
U$, $L_\lambda$ is admissible. We thus assume this holds for all $\lambda\in\Lambda$.

An element in the tangent space of $T_\lambda U$ is simply a vector field
 $\tilde Y_\lambda$ along
$L_\lambda$.
We express the corresponding vector field along
$\Pi_P(L_\lambda)$, $\Pi_P(\tilde Y_\lambda),$ as $Y_\lambda$ (or sometimes
$Y_\lambda^\sigma$ to encode more information mentioned below).
Consider the bundles
\begin{eqnarray*}
\cand_{2,\epsilon,\Lambda}({\mathbf
b}) & = &
\bigcup_{\lambda\in\Lambda}\cand_{2,\epsilon}({\mathbf b}),
\\
\E_\Lambda& =& \bigcup_{\lambda\in\Lambda}\E(L_\lambda)
\end{eqnarray*}
To find local trivializations of these bundles we proceed as follows. Denote the Legendrian
embeddings corresponding to $L_\lambda$ by
$$
k_\lambda\colon L\to P\times \R,
$$
where $L=L_0$. Let $c_j^\pm[\lambda]\subset L$ be the Reeb chord end points of $L_\lambda$ and write
$c_j^\pm=c_j^\pm[0]$. Let $\phi_\lambda\colon L\to L$ be diffeomorphisms such that
$$
\phi_\lambda(c_j^\pm)=c_j^\pm[\lambda],
$$
for all $c_j^\pm[\lambda]$, which is the identity outside $B(c_j^\pm, 9\delta)$. Construct families
$\Phi_\lambda^\sigma$ of diffeomorphisms of $P$ such that the following holds
\begin{itemize}
\item $\Phi^\sigma_\lambda(c_j^\ast[0])= c_j^\ast[\lambda]$ for all
$\sigma\in[0,1]$.
\item For all $x\in L\setminus(\bigcup B(c_j^\pm;10\delta))$
$$
J \circ d\Phi^\sigma_\lambda = d\Phi^\sigma_\lambda\circ J,
$$
at $\Pi_P(k_\lambda(x))$ for all $\sigma\in[0,1]$.
\item For $x\in L\setminus B(c_j^-,5\delta)$,
$$
\Phi^0_\lambda(\Pi_P(x))=\Pi_P(k_\lambda(x)).
$$
\item For $x\in L\setminus B(c_j^+,5\delta)$,
$$
\Phi^1_\lambda(\Pi_P(x))=\Pi_P(k_\lambda(x)).
$$
\item For $x\in L\setminus B(c_j^-,7 \delta)$,
$$
J \circ d\Phi^0_\lambda = d\Phi^0_\lambda \circ J,
$$
at $\Pi_P(k_\lambda(x))$.
\item For $x\in L\setminus B(c_j^+,7 \delta)$,
$$
J \circ d\Phi^1_\lambda = d\Phi^1_\lambda \circ J,
$$
at $\Pi_P(k_\lambda(x))$.
\end{itemize}

To define local coordinates on the manifold $\cand_{2,\epsilon,\Lambda}({\mathbf b})$ we first note
that we can cover $\cand_{2,\epsilon,0}({\mathbf b})$ by charts centered at maps which are constant
close to their punctures. If $(u,f,\kappa)$ is an element with $u$ and $f$ constant close to its
punctures and $U$ is the corresponding local coordinates define the map
$$
\Psi\colon U\times\Lambda\to\cand_{2,\epsilon, \Lambda}({\mathbf b}),
$$
$$
\Psi(v,\lambda)(\zeta)=
\exp_{\Phi_\lambda^{\rho(F(\zeta))}(u(\zeta))}^{\rho(F(\zeta))}(d\Phi^{\rho(F(\zeta))}_\lambda
v(\zeta)),\quad \zeta\in D_m.
$$
\begin{lma} The map $\Psi(v,\lambda)$ lies in
$\cand_{2,\epsilon,\Lambda}({\mathbf b})$, over $\lambda\in\Lambda$.
\end{lma}
\begin{pf}
The fact that the $(2,\epsilon)$-norm of $\Psi(v,\lambda)$ is bounded is straightforward. We
check that $\Psi(v,\lambda)$ satisfies the boundary conditions. We note first that $\Psi(0,\lambda)$
does: for $\zeta\in\pa D_m$ and $v$ sufficiently small either $\rho(F(\zeta))$ is constantly equal
to $1$ or $0$ or $\Phi^{\rho(F(\zeta))}_\lambda\circ k_\lambda$ is independent of $\zeta$. Thus, at
such $\zeta$
$$
d\Phi^\sigma_\lambda\circ du + J d\Phi^\sigma_\lambda\circ du\circ j_\kappa=
d\Phi^\sigma_\lambda\circ (du + J\circ du\circ j_\kappa)=0,
$$
by the above properties of $\Phi^\sigma_\lambda$. The lemma then follows from Lemma~\ref{lma.coords} once we
establish $\bar\nabla_J d\Phi^\sigma_\lambda v=0$ on the boundary. We compute in complex coordinates
$\tau+ i t$ where $\tau$ is tangent to the boundary
\begin{align*}
\nabla_\tau (d\Phi^\sigma_\lambda v) + J\nabla_t (d\Phi^\sigma_\lambda v)j_\kappa =
\nabla_\tau X + J\nabla_t X,
\end{align*}
where $X$ is tangent to $L_\lambda$. Since $\Psi(0,\lambda)$ is holomorphic on the boundary we find
that
$$
\pa_t=-J\pa_\tau,
$$
where $\pa_\tau$ and $\pa_t$ denotes the derivatives of $\Psi(0,\lambda)$ in directions parallel to
the boundary and orthogonal to the boundary respectively. Moreover, $\pa_\tau$ is tangent to $L_\lambda$. Thus the boundary conditions hold since
$$
\nabla_{JT} X= J\nabla_T X,
$$
where $T$ is any vector field tangent to $L$ along $L$ and where $\nabla$ is the Levi-Civita connection of the metric $\hat g$ on $TL$, see Lemma 4.6 in \cite{ees2}.
\end{pf}

Finally we must trivialize the bundle $\E_\Lambda$ over such a coordinate region. To this end we consider
the auxiliary complex structures on $u^\ast TP$
$$
J_{\lambda,v}=
[\Pi_{d\Phi_\lambda^\sigma(v)}]^{-1}[d\Phi_\lambda^\sigma]^{-1}
J
[\Pi_{d\Phi_\lambda^\sigma(v)}][d\Phi^\sigma],
$$
which corresponds to the complex structure $J$ over $\Psi(v,\lambda).$
We may now proceed in the
standard way to trivialize: writing
$$
J_{\lambda,v}=J(1+S_{\lambda,v})(1-S_{\lambda,v})^{-1},
$$
and mapping $(j_\kappa,J)$-complex anti-linear maps $A$ to $(j_\kappa, J_{\lambda,v})$-complex anti-linear maps $(1-S_{\lambda,v})A$.

\subsection{Varying the almost complex structure}\label{DJ.sec}

Let $L$ be a fixed Legendrian submanifold and let $J_\lambda$ be a family of almost complex structures
on $P$ compatible with $\omega$ and adapted to $L$ and such that $L$ is admissible with respect to
$J_\lambda$ for each $\lambda\in \Lambda$. Consider the bundles
\begin{eqnarray*}
\cand_{2,\epsilon,\Lambda}({\mathbf b})& =&\bigcup_{\lambda \in \Lambda} \cand_{2,\epsilon}({\mathbf b},J_\lambda)\\
\E_\Lambda& = &\bigcup_{\lambda\in\Lambda}\E(J_\lambda)
\end{eqnarray*}
Note we recycle notation from Section \ref{DLs.sec}, relying on the context
to distinguish bundles for deforming Legendrian submanifolds versus deforming
almost complex structures.
We first consider local coordinates on $\cand_{2,\epsilon, \Lambda}({\mathbf b})$. Recall that the metric
$g^\sigma$ was constructed using an immersion $S\colon U \to P,$ of a neighborhood $U$ of
the zero section in $TL,$ which is holomorphic along the $0$-section with respect to the canonical complex structure along the $0$-section and the almost complex structure $J$ on $P$. Hence it depends on $J$ and as $J$ changes we must change this immersion.
Note that we assume that $J_\lambda$ is fixed in some neighborhood of each double point of the
Lagrangian immersion. We now fix a family of $1$-parameter families of diffeomorphisms
$$
\Phi_\lambda\colon P\to P,
$$
such that
\begin{itemize}
\item $\Phi_\lambda$ fixes a small neighborhood $W_j$ of each $c_j^\ast$ where $J$ does not vary.
\item Let $U$ be the intersection of a regular neighborhood of $L$ with the complement of the
neighborhoods fixed by $\Phi_\lambda$. This is a tubular neighborhood of
$\Pi_P(L)\setminus(\cup_j W_j)$. Let $\Phi_\lambda$ satisfy $\Phi_\lambda= S_\lambda\circ S_0^{-1}$ in a
smaller neighborhood contained in this neighborhood and $\Phi_\lambda=\id$ outside the neighborhood.
\end{itemize}

We can now take the local coordinates as
$$
(u,\lambda)\mapsto \Phi_\lambda\circ u.
$$
On the bundle $\E_\Lambda$ we use $d\Phi_\lambda$ for conjugation and proceed as before.

\begin{rmk}\label{commute.rmk}
The technical reason for separating the variations of the almost complex structure from the deformations of the Legendrian submanifold in time, see Lemma \ref{lmaadm1}, is the following. In order to extend the bundle of $\E$ to a bundle $\E_\Lambda$ over a space $\Lambda=\Lambda_L\times\Lambda_J$ which parameterizes simultaneous deformations of the Legendrian submanifold $L$ and the almost complex structure $J$ we must define suitable diffeomorphisms $\Phi_\lambda\colon P\to P$ which are holomorphic along the image of $L$ under $\Pi_P$ as described above in the cases when $L$ and $J$ vary separately. If $\frac{\pa}{\pa \lambda_L}$ and $\frac{\pa}{\pa\lambda_J}$ are vector fields at $(0,0)\in\Lambda$ in $\Lambda_L$- and $\Lambda_J$-directions respectively the vector fields $\frac{\pa}{\pa\lambda_L}\Phi_\lambda|_{\lambda=0}$ and $\frac{\pa}{\pa\lambda_J}\Phi_\lambda|_{\lambda=0}$ in a neighborhood of $\Pi_P(L)\subset P$ do, in general, not commute. Therefore it is not straightforward to produce local trivializations for the bundle $\E_\Lambda$ for simultaneous and independent deformations of $L$ and $J$.
\end{rmk}

%----------------------------------------------------------------------
\subsection{The equation for holomorphic disks and its linearization}\label{linearizeJ.sec}
%----------------------------------------------------------------------
In this section we compute the linearization of the $\bar\pa_J$-equation in the coordinates
described above.

\subsubsection{Linearization when $L$ varies} Consider the $\bar\pa_{J,j}$-operator as giving a
section of the bundle
$$
\E_\Lambda\to\cand_{2,\epsilon,\Lambda}({\mathbf b}).
$$
To compute the linearization we use the trivialization described in Section~\ref{sec.coords}: in
local coordinates $U\times\Lambda \subset u^\ast(TP)\times\Lambda$ we trivialized the bundle
$\E_\Lambda$ by identifying all nearby fibers with
$$
\sblv_{1,\epsilon}(D_m, \Hom^{0,1}(TD_m,u^\ast(TP)).
$$
To simplify notation and make it consistent with Section \ref{DLs.sec},
we write $Y^\sigma_\lambda=\frac{\pa}{\pa\lambda}\Phi^\sigma_\lambda$, where
$\sigma\in[0,1]$ is as in Section \ref{DLs.sec} and
$\lambda\in\Lambda$ is a coordinate on the space of Legendrian submanifolds. Thus $Y_\lambda^\sigma$ is a $\sigma$-dependent vector filed generating the $\sigma$-dependent diffeomorphisms $\Phi^\sigma_\lambda$ along a curve in $\Lambda$ with tangent vector $\frac{\pa}{\pa\lambda}$.
We write $\frac{d}{d\sigma} Y^\sigma_\lambda$ for the vector field
$\lim_{\tau\to\sigma}\frac{1}{\tau-\sigma} (Y^\sigma_\lambda-Y^\tau_\lambda)$ and we let $\Lie_V$ denote
the Lie derivative with respect to the vector field $V$. As in \ref{confstr.sssec}, we consider $T_\kappa\conf_m$ as a subspace of $\End^{0,1}(TD_m)$ and we write $L_\lambda$ for the Legendrian submanifold corresponding to $\lambda\in\Lambda$.

Consider the boundary condition
\begin{equation}\label{5.36.eq}
\int_{\pa D_m}\la v,w\ra\,ds=0,  \text{ for all $w\in C^0_0(\pa D_m,T^{0,1}D_m\otimes\C^n)$.}
\end{equation}
and define for $k \ge 1$
\begin{equation}
\sblv_{k,\mu}[0](D_m,u^\ast(TP))  =  
\bigl\{v\in\sblv_{k, \mu}(D_m, u^\ast(TP)) \,\colon v \text{ satisfies \eqref{5.36.eq}} \bigr\} 
 \label{H1mu.eq} 
\end{equation}
\begin{equation}
\sblv_{k,\mu}[0]({T^\ast}^{0,1}D_m\otimes\C^n)  = 
\bigl\{v\in\sblv_{k, \mu}(D_m, \text{Hom}^{0,1}
(TD_m, u^\ast (TP))) \,\colon v \text{ satisfies \eqref{5.36.eq}} \bigr\}
 \notag
\end{equation}

\begin{lma}\label{linearization1.lem}
Let
$$
(v,\lambda,\gamma)\in \sblv_{2,\epsilon}[0](D_m,u^\ast(TP))\times T_0\Lambda\times \End^{0,1}(TD_m)
$$
be a tangent vector of $\cand_{2,\epsilon,\Lambda}({\mathbf b})$ at $(u,f,\kappa;L_0)$, let $J$ be the
almost complex structure on $P$, and let $j$ be the complex structure on $D_m$ corresponding to
$\kappa$. Then
\begin{align}\notag
D\,\bar\pa_{J,j}[v,\lambda,\kappa]=&
\Bigl(\bar\nabla_{J,j}\,v
-\tfrac12\,J\circ (\nabla_v J)\circ\pa_{J,j}\,u\Bigr)\\\notag
&+\Bigl(
-\tfrac12J\circ (\Lie_{Y^\sigma_\lambda} J)\circ \pa_{J,j}\, u
+ \tfrac{d}{d\sigma}Y_\lambda^\sigma\cdot d\sigma +
J\tfrac{d}{d\sigma}Y_\lambda^\sigma\cdot d\sigma\circ j
\Bigr)\\\label{eqlin1}
&-\pa_{J,j}\, u\circ \gamma.
\end{align}
\end{lma}

\begin{rmk}
The fact that the local coordinates and the $\bar\pa_{J,j}$-section are $C^1$ can be proved in a way
similar to the corresponding results in \cite{ees2}.

In the case of varying Legendrian submanifold we added the condition that the center $u$ of a coordinate patch
be a map constant near the punctures. Note however that \eqref{eqlin1} makes sense for all elements in
$\cand_{2,\epsilon,\Lambda}({\mathbf b})$. The $C^1$-properties of the local coordinates together
with an approximation argument show that \eqref{eqlin1} holds in general.
\end{rmk}

\begin{pf}
Using the fact that the section and the local coordinates are $C^1$ it is sufficient to compute the
partial derivatives in order to compute the differential.

Consider first the partial derivative with respect to $v$. Letting $\exp^\sigma$ denote the
exponential map in the metric $g^\sigma$ we must compute
\begin{equation}\label{eqfullv}
\bar\pa_{J,j}\exp^{\sigma(\zeta)}_{u(\zeta)}(\epsilon v(\zeta)),
\end{equation}
to first order in $\epsilon$ in our trivialization.
(This "$\epsilon$" is unrelated to the "$\epsilon$" from the weight functions,
e.g., $\cand_{2,\epsilon,\Lambda}$.)
First note that if $x\in P$ and $\xi\in T_xP$
then
$$
\frac{d}{d\sigma}\exp_x^\sigma(\epsilon \xi)=\Ordo(\epsilon^2).
$$
Therefore we may view $\sigma$ as fixed in our computation below. In our chosen trivialization \eqref{eqfullv} is expressed as
\begin{equation}\label{eqpartv}
(1-S_{\epsilon v})^{-1}\,\Pi_{-\epsilon v}\,\Bigl(d\exp^\sigma_{u(\zeta)}(\epsilon v(\zeta))+J\circ
d\exp^\sigma_{u(\zeta)}(\epsilon v(\zeta))\circ j\Bigr),
\end{equation}
where $S_{\epsilon v}=(J+J_{\epsilon v})^{-1}(J_{\epsilon v}-J)$. Let $\eta\in T_{\zeta} D_m$ then
$$
d[\exp^\sigma_{u}(\epsilon v)]\cdot\eta=X(\epsilon),
$$
where $X$ is the Jacobi field along $\gamma(s)=\exp_{u(\zeta)}^\sigma(sv(\zeta))$ with
$X(0)=du\cdot\eta$ and $\nabla_s X(0)=\nabla_{du\cdot\eta}\, v$. If $E_1,\dots, E_{2n}$ is a
parallel frame along $\gamma$ and we write $X=\Sigma_j e^jE_j$ then the Jacobi equation is
$$
\frac{d^2 e^k}{ds^2}+ \Sigma_j e^j\rho_j^k=0,
$$
where $R(E_j,\dot\gamma)\dot\gamma=\Sigma_k \rho_j^k E_k$. Thus,
\begin{equation}\label{eqa}
\Pi_{-\epsilon v}\Bigl(
d[\exp^\sigma_{u}(\epsilon v)]\cdot\eta\Bigr)=
du\cdot\eta+\epsilon \nabla_{du\cdot\eta} v +\Ordo(\epsilon^2).
\end{equation}
A similar calculation, using also the fact that $\Pi_{-\epsilon v}J-J\Pi_{-\epsilon v}=\epsilon \nabla_v J+\Ordo(\epsilon^2)$ gives
\begin{align}\notag
\Pi_{-\epsilon v}\Bigl(J \left(d[\exp^\sigma_{u}(\epsilon
v)]\cdot j\eta\right)\Bigr) =& J \left(du\cdot j\eta\right)\\
\label{eqaa}
&+\epsilon \Bigl(J\left(\nabla_{du\cdot j\eta} v\right) +
\nabla_v J\left(du\cdot j\eta\right)\Bigr) \\ & \notag +\Ordo(\epsilon^2).
\end{align}
We next we determine $(1-S_{\epsilon v})^{-1}$ up to first order in $\epsilon$.
Note that $(1-S_{\epsilon v})^{-1}=1+S_{\epsilon v}+\Ordo(\epsilon^2)$ and that
$J_{\epsilon v}=J+\Ordo(\epsilon)$. Thus, the definition of $S_{\epsilon v}$ implies
\begin{align}\notag
S_{\epsilon v} &=(J + J + \Ordo(\epsilon))^{-1}
(\Pi_{-\epsilon v}\circ J\circ \Pi_{\epsilon v} - J) \\
\notag
&=-\frac12 J(\Pi_{-\epsilon v}\circ J\circ
\Pi_{\epsilon v} - J)+\Ordo(\epsilon^2)\\
\label{eqaaa}
&=-\epsilon\frac12 J\circ \nabla_vJ +\Ordo(\epsilon^2).
\end{align}
The partial derivative of \eqref{eqfullv} with respect to $v$ is the first order term in the
$\epsilon$-expansion of \eqref{eqpartv}. Using \eqref{eqa}, \eqref{eqaa}, and \eqref{eqaaa}, we find
that
\begin{align*}
D\,\bar\pa_{J,j}[(v,0,0)]&=
\bar \nabla_{J,j}\, v + (\nabla_vJ)(du)j -\frac12 J(\nabla_vJ)du
- \frac12 J(\nabla_vJ)J(du)j\\
&=\bar \nabla_{J,j}\, v +\frac12(\nabla_vJ)\bigl((du)j + J(du)\bigr)\\
&=\bar \nabla_{J,j}\, v -\frac12 J(\nabla_vJ)(\pa_{J,j}\,u),
\end{align*}
which agrees with the first term in \eqref{eqlin1}.

Consider next the partial derivative with respect to $\lambda$. We must compute
\begin{equation}\label{eqfulllambda}
\bar\pa_{J,j}\Phi_{\epsilon\lambda}^\sigma(u)
\end{equation}
to first order in $\epsilon$ in our trivialization where \eqref{eqfulllambda} is expressed as
\begin{equation}\label{eqpartlambda}
(1-S_{\epsilon\lambda, 0})^{-1}\left(D\Phi_{-\epsilon\lambda}^\sigma\right)
\left(d\Phi_{\epsilon\lambda}^\sigma(u)+
J\circ d\Phi_{\epsilon\lambda}^\sigma(u)\circ j\right),
\end{equation}
where $D\Phi_{-\epsilon\lambda}^\sigma$ denotes the differential of the diffeomorphism
$\Phi_{-\epsilon\lambda}^\sigma\colon P\to P$.

As above we observe that $(1-S_{\epsilon\lambda, 0})^{-1}=1+S_{\epsilon\lambda, 0}+\Ordo(\epsilon^2)$.
We start by computing $S_{\epsilon\lambda, 0}.$ By definition
\begin{align}\notag
S_{\epsilon\lambda, 0} &=-\frac12 J(D\Phi^\sigma_{-\epsilon\lambda}\circ J\circ
D\Phi^\sigma_{\epsilon\lambda}-J)+\Ordo(\epsilon^2)\\\label{eqb}
&=-\epsilon\frac12 J\circ \left(\Lie_{Y^\sigma_\lambda}J\right) +\Ordo(\epsilon^2).
\end{align}
Next, in local coordinates $x\in U\subset\R^{2n}$ on $P$,
$$
\Phi_{\epsilon\lambda}^\sigma(x)=x+
\epsilon\, Y_\lambda^\sigma(x)+\Ordo(\epsilon^2).
$$
Therefore,
$$
D\Phi_{\epsilon\lambda}^\sigma=\id+\epsilon\,
D Y_\lambda^\sigma+\Ordo(\epsilon^2),
$$
and thus
\begin{equation}\label{eqbb}
d(\Phi_{\epsilon\lambda}^\sigma(u))=D\Phi_{\epsilon\lambda}^\sigma\cdot
du+ \epsilon \left(\tfrac{d}{d\sigma}Y_\lambda^\sigma\right)\cdot
d\sigma+\Ordo(\epsilon^2)
\end{equation}

The partial derivative with respect to $\lambda$ is the first order (in $\epsilon$) term in
\eqref{eqpartlambda}. Using \eqref{eqb} and \eqref{eqbb} we find that
\begin{align*}
D\,\bar\pa_{J,j}[(0,\lambda,0)]=&
-\frac12 J(\Lie_{Y^\sigma_\lambda}J)(du + J(du)j)
+(\Lie_{Y_\lambda^\sigma}J)(du)j\\
&+(\tfrac{d}{d\sigma}Y_\lambda^\sigma)d\sigma
+J(\tfrac{d}{d\sigma}Y_\lambda^\sigma)(d\sigma)j\\
=&-\frac12J(\Lie_{Y_\lambda^\sigma}J)(du-J(du)j)
+(\tfrac{d}{d\sigma}Y_\lambda^\sigma)d\sigma
+J(\tfrac{d}{d\sigma}Y_\lambda^\sigma)(d\sigma)j
\end{align*}
which agrees with the second term in \eqref{eqlin1}.

Finally, we consider the partial derivative with respect to $\gamma$. It equals the order $\epsilon$
term in
\begin{equation}\label{eqfullgamma}
\bar\pa_{J,j_{\epsilon\gamma}}\,u,
\end{equation}
where $\epsilon\gamma$ denotes the complex structure
$j_{\epsilon\gamma}=j_\kappa(1+\epsilon\gamma)(1-\epsilon\gamma)^{-1}$ on $D_m$. In our
trivialization \eqref{eqfullgamma} is expressed as
$$
\Bigl(du + J\circ du \circ
j(1+\epsilon\gamma)(1-\epsilon\gamma)^{-1}\Bigr)(1+\epsilon\gamma)^{-1}.
$$
Thus, arguing as above, we find
$$
D\,\bar\pa_{J,j}[(0,0,\gamma)]=-\pa_{J,j}\,u\circ\gamma,
$$
which agrees with the third term in \eqref{eqlin1}. This finishes the proof of the lemma.
\end{pf}

\subsubsection{Linearization when $J$ varies}
Consider the $\bar\pa_{J,j}$-operator as giving a section of the bundle
$$
\E_\Lambda \to\cand_{2,\epsilon,\Lambda}.
$$
To compute the linearization we use the trivialization described above: in local coordinates
$U\times\Lambda \subset u^\ast(TP)\times\Lambda$ we trivialized the bundle $\E_\Lambda$ by identifying
all fibers with
$$
\sblv_{1,\epsilon}(D_m, \Hom^{0,1}(TD_m,u^\ast(TP))).
$$
To simplify notation and make it analogous to Section \ref{DLs.sec} we write $Y_\lambda=\frac{\pa}{\pa\lambda}\Phi_{\lambda}$, where
$\lambda\in\Lambda$ is the coordinate on the
space of adapted almost complex structures. Thus $Y_\lambda$ is a vector field generating the diffeomorphisms $\Phi_\lambda$ along a curve in $\Lambda$ with tangent vector $\frac{\pa}{\pa\lambda}$. Let $\Lie_V$
denote the Lie derivative with respect to the vector field $V$ and let
$J_\lambda=J_0(1+S_\lambda)(1-S_\lambda)^{-1}$, where $S_\lambda$ anti-commutes with $J_0$.
\begin{lma}\label{linearization2.lem}
Let
$$
(v,\lambda,\gamma)\in \sblv_{2,\epsilon}[0](D_m,u^\ast(TP))\times T_0\Lambda\times \End^{0,1}(TD_m)
$$
be a tangent vector of $\cand_{2,\epsilon,\Lambda}({\mathbf b})$ at $(u,f,j;J_0)$. Then
\begin{align}\notag
D\,\bar\pa_{J_0,j}[v,\lambda,\kappa]=&
\Bigl(\bar\nabla_{J,j}\,v
-\tfrac12\,J_0\circ (\nabla_v J_0)\circ\pa_{J_0,j}\,u\Bigr)\\\notag
&+\Bigl(
-(\tfrac12J_0\circ (\Lie_{Y_\lambda} J_0)+S_\lambda)\circ \pa_{J_0,j}\, u
\Bigr)\\\label{eqlin2}
&-\pa_{J,j}\, u\circ \gamma.
\end{align}
\end{lma}
\begin{pf}
The proof is analogous to the proof of Lemma~\ref{linearization1.lem}. In fact the appearance of the
first and last terms in \eqref{eqlin2} is a consequence of the proof of
Lemma~\ref{linearization1.lem}. To finish the proof we must compute the partial derivative with
respect to $\lambda$. In other words, we must compute
\begin{equation}\label{eqfullJ}
\bar\pa_{J_{\epsilon\lambda},j}\Phi_{\epsilon\lambda}(u)
\end{equation}
to first order in $\epsilon$ in our trivialization, where \eqref{eqfullJ} is expressed as
$$
(1-S_{\epsilon \lambda})^{-1}(D\Phi_{-\epsilon\lambda})(d\Phi_{\epsilon \lambda}(u)+J_{\epsilon \lambda}(d\Phi_{\epsilon \lambda}(u))j).
$$
First note that
$$
J_{\epsilon\lambda}=J_0(1+S_{\epsilon\lambda})(1-S_{\epsilon\lambda})^{-1}
=J_0(1+2 S_{\epsilon \lambda}+\Ordo(\epsilon^2)).
$$
We next compute $S_{\epsilon \lambda}$. We have
\begin{align*}
S_{\epsilon \lambda} &=-\frac12 J_0(d\Phi_{-\epsilon\lambda}\circ J_{\epsilon \lambda}\circ
d\Phi_{\epsilon \lambda} - J_0)+\Ordo(\epsilon^2)\\
&=\epsilon J_0\circ\left(-\frac12\Lie_YJ_0 \right)+\Ordo(\epsilon^2).
\end{align*}

We thus find
\begin{align*}
D\,\bar\pa_{J_0,j}[0,\lambda,0]=&
-\tfrac12J_0(\Lie_{Y_\lambda}J_0)(du)
-\tfrac12J_0(\Lie_{Y_\lambda}J_0)J_0(du)j\\
&+\Lie_{Y_\lambda}J_0(du)j + 2J_0S_\lambda(du)j\\
&+S_\lambda(du)+S_\lambda J_0(du)j\\
=&-\tfrac12 J_0(\Lie_{Y_\lambda}J_0)(\pa_{J_0,j}\, u)
+S(\pa_{J_0,j}\, u),
\end{align*}
which agrees with the second term in \eqref{eqlin2}. The lemma follows.
\end{pf}

\begin{rmk}
Consider Equations \eqref{eqlin1} and \eqref{eqlin2}. Assume that we are linearizing at a point in $\cand_{2,\epsilon,\Lambda}(b)$ where the almost complex structure in $P$ is $J$ and that in $D_m$ is $j$. Then the right hand sides of \eqref{eqlin1} and \eqref{eqlin2} are $(J,j)$-complex anti-linear.
Moreover their restrictions to the boundary vanishes by our choices of $\Phi_\lambda^\sigma$ and $\Phi_\lambda$, respectively.
\end{rmk}

%----------------------------------------------------------------
\section{Manifolds and gluing}\label{mainann.sec}
%----------------------------------------------------------------
In this section, we show how to adopt the techniques of Sections 5, 6 and 7 of \cite{ees2} and the orientations in \cite{ees3}, in the setup described in Section \ref{ban.sec}, to prove Proposition~\ref{mfld.prop} and  Lemma~\ref{dsquared=0.lem}. In particular, we discuss Fredholm  properties in Subsection \ref{fred.sec}, transversality in Subsection \ref{trans.sec}, gluing results in Subsection \ref{glue.sec}, and orientations on the moduli spaces in Subsection \ref{spin.sec}.

%----------------------------------------------------------------
\subsection{Fredholm properties} \label{fred.sec}
%----------------------------------------------------------------
Let $(u,f) \in \cand_{2,\epsilon}({\mathbf b}; \kappa, J)$ be a map with boundary conditions at some
admissible $L,$ that is, $(u,f)$ satisfies (\ref{boundary1.eq}) and (\ref{boundary2.eq}). We wish to
show that the linearization of $\bar{\pa}_J$ at $u$ is Fredholm and of a certain index.

We choose a complex trivialization of the pull back bundle $u^\ast(TP)\rightarrow D_m.$ In this
trivialization, let $A: \pa D_m \rightarrow \U(n)$ represent $u^\ast |_{\pa D_m}(TL).$ (The Lagrangian subspace associated to $A$ is $A\cdot\R^n\subset\C^n$.) Recall that near the image under $u$ of each puncture our set-up looks like $\C^n$; thus, like in \cite{ees2}, $A$ is
{\em{small at infinity}} and over each puncture $p_j$, $A^+_j \R^n$ and $A^-_j \R^n$ are transverse.
Here $A^\pm_j$ denotes the two limiting matrices near the $j$-th puncture and ``small at infinity''
means $A$ converges in the $C^2$-norm to $A^\pm_j$ as $z \in \pa D_m$ converges to $p_j.$ See
Definition 5.7 \cite{ees2} for more details.

As in Section 5 of \cite{ees2}, we consider the stationary case only (no changes in $J$ or $L$) and
we ignore the last term $-\pa_{J,j}\, u\circ \gamma$ in (\ref{eqlin1}) and (\ref{eqlin2})
associated to varying the conformal structure. These variations are finite dimensional and effects
the Fredholm index in a standard way. In this case, (\ref{eqlin1}) and
(\ref{eqlin2}) reduce to
\begin{equation}
\label{lin.eq}
D\,\bar\pa_{J,j}[v]= \bar\nabla_{J,j}\,v -\tfrac12\,J\circ (\nabla_v J)\circ\pa_{J,j}\,u :=
\bar\nabla_{J,j}\,v + Kv
\end{equation}
where $v \in \sblv_{2,\epsilon}[0](D_m,u^\ast(T^\ast P)) $ is a tangent vector of
$\cand_{2,\epsilon}({\mathbf b})$ at $(u,f).$

\begin{lem}
\label{K.lem}
The operator $K$ vanishes on $\pa D_m$.
\end{lem}

\begin{proof}
Let $T$ be any vector tangent to $L$ and $X$ any vector field in $T(TL).$ Then
\[
(\nabla_TJ)X = \nabla_T(JX) - J(\nabla_TX) = 0,
\]
by (4.13) in the proof of Lemma 4.6 in \cite{ees2}.
\end{proof}

Consider the boundary condition
\begin{align}\notag
&\int_{\pa D_m}\la v,w\ra\,ds=0,
\text{for all $w\in C^0_0(\pa D_m,\C^n)$}  \\  \label{5.35.eq}
& \text{such that
$w(\zeta)\in iA(\zeta)\R^n$ for all $\zeta \in\pa D_m$}
\end{align}
and define
\begin{align}\notag
\sblv_{2,\mu}(\C^n;A)=
\bigl\{v\in&\sblv_{2,\mu}(\C^n)\colon v\text{ satisfies \eqref{5.35.eq}, }
\bar\nabla_{J,j} v + Kv,\\\label{H2mu.eq}
&\text{ satisfies \eqref{5.36.eq} }\bigr\}.
\end{align}

By Lemma~\ref{K.lem}, the definition of $\sblv_{2,\mu}(\C^n;A)$ does not change if we were to
replace ``$\nabla_{J,j} v + Kv $ satisfies \eqref{5.36.eq}'' with ``$\nabla_{J,j} v $ satisfies
\eqref{5.36.eq},'' and so we will alternate between the two conditions.

\begin{lem}\label{cptK.lem}
$K : \sblv_{2,\mu}(\C^n;A) \rightarrow \sblv_{1,\mu}[0]({T^\ast}^{0,1}D_m\otimes\C^n)$
is a compact operator.
\end{lem}
\begin{proof}
In local coordinates around a point in $P$ which lies in the image of $u$ we may think of $Kv$ as a matrix valued function. Its components are finite sums of terms of the form $h(u(\zeta))v_j(\zeta)(\pa_r u_k)(\zeta)$, where $h$ is a smooth function with bounded derivatives of all orders which depends on the metric $g^\sigma$, where $v_j$ and $u_k$ denote components of $v$ and $u$ in the given coordinate system, and where $\pa_r$, $r=1,2$, denotes partial derivative with respect to the variables in the disk. Using a partition of unity argument together with the fact that, by definition of $\cand_{2}({\mathbf b})$, $u$ takes a neighborhood of each puncture in $D_m$ into a neighborhood of a double point of $\Pi_P(L)$, the lemma follows once we establish the following claim.

\begin{clm}
Let ${\bf M}(2n,2)^\ast$ denote the space of linear forms on the linear space of $(2n\times 2)$-matrices and let $(\R^{2n})^\ast$ denote the dual space of $\R^{2n}$. Let $h\colon \R^{2n}\to {\bf M}(2n,2)^\ast\otimes(\R^{2n})^\ast$ be a smooth function with all derivatives bounded. Fix $u\in \sblv_{2}(D_m,\R^{2n})$.
Then
$$
F\colon\sblv_{2}(D_m,\R^{2n})\to \sblv_{1}(D_m,\R),\quad
F(v)=h(u)\cdot du\cdot v,
$$
where $\cdot$ denotes contraction of tensors, is a compact operator.
\end{clm}

Let $w\colon D_m\to\R^{2n}$ be a smooth function with compact support. Define the auxiliary operator
$$
F_w\colon\sblv_{2}(D_m,\R^{2n})\to\sblv_1(D_m,\R),\quad
F_w(v)=h(u)\cdot dw\cdot v.
$$
To establish the claim we first derive estimates for the operator norm of the difference $F-F_w$. To this end we note that the $(2,2)$-norm on $D_m$ controls the supremum norm,
$\|v\|_\infty\le C\|v\|_2$,
and that the $(1,2)$-norm controls the $(0,p)$-norm for all $2\le p<\infty$.
We calculate
\begin{align}\notag
&\int_{D_m}|F(v)-F_w(v)|^2\,dA\le\int_{D_m}|h(u)|^2|du-dw|^2|v^2|\,dA\le\\\label{eq0Dnorm}
&C\left(\int_{D_m}|du-dw|^4\,dA\right)^{\tfrac12}\left(\int_{D_m}|v|^4\,dA\right)^{\tfrac12}
\le C\|u-w\|_2^2\|v\|_2^2.
\end{align}
(In this proof, for simplicity of notation, we will let $C$ denote a constant which is allowed to change its value but which is always independent of $u$ and $v$.) We derive a similar bound for the $L^2$-norm of the derivative:
\begin{align*}
d\bigl(F(v)-F_w(v)\bigr) &= Dh(u)\cdot du\cdot (du-dw)\cdot v\\
&+h(u)\cdot(d^2u-d^2w)\cdot v+h(u)\cdot(du-dw)\cdot dv.
\end{align*}
Therefore
\begin{align*}
\int_{D_m}&|d(F(v)-F_w(v))|^2\,dA\le\\
&C\left(\int_{D_m}|du|^2|du-dw|^2|v|^2\,dA\right.
+\int_{D_m}|d^2u-d^2w|^2|v|^2\,dA\\
&\quad+ \int_{D_m}|du-dw|^2|dv|^2\,dA
+\int_{D_m}|du||du-dw||d^2u-d^2w||v|^2\,dA\\
&\quad+\int_{D_m}|du||du-dw|^2|v||dv|\,dA
\left.+\int_{D_m}|du-dw||d^2u-d^2w||v||dv|\,dA\right).
\end{align*}
We estimate the terms on the right hand side one by one.
\begin{align}\notag
\int_{D_m}|du|^2|du-dw|^2|v|^2\,dA&\le
\|v\|_\infty^2\left(\int_{D_m}|du|^4\,dA\right)^{\tfrac12}
\left(\int_{D_m}|du-dw|^4\,dA\right)^{\tfrac12}\\\label{eq1DnormI}
&\le C\|u-w\|_2^2\|u\|_2^2\|v\|_2^2,
\end{align}

\begin{align}\label{eq1DnormII}
\int_{D_m}|d^2u-d^2w|^2|v|^2\,dA\le
\|v\|_\infty^2\|u-w\|_2^2\le C\|u-w\|_2^2\|v\|_2^2,
\end{align}

\begin{align}\notag
\int_{D_m}|du-dw|^2|dv|^2\,dA&\le
\left(\int_{D_m}|du-dw|^4\,dA\right)^{\tfrac12}
\left(\int_{D_m}|dv|^4\,dA\right)^{\tfrac12}\\\label{eq1DnormIII}
&\le C\|u-w\|_2^2\|v\|_2^2,
\end{align}

\begin{align}\notag
&\int_{D_m}|du||du-dw||d^2u-d^2w||v|^2\,dA\le\\\notag
&\|v\|_\infty^2\left(\int_{D_m}|du|^2|du-dw|^2\,dA\right)^{\tfrac12}
\left(\int_{D_m}|d^2u-d^2w|^2\,dA\right)^{\tfrac12}\\\label{eq1DnormIV}
&\le C\|u-w\|_2^2\|u\|_2\|v\|_2^2,
\end{align}

\begin{align}\notag
\int_{D_m}|du||du-dw|^2|v||dv|\,dA&\le
\|v\|_\infty\left(\int_{D_m}|du|^2|du-dw|^4\,dA\right)^{\tfrac12}
\left(\int_{D_m}|dv|^2\,dA\right)^{\tfrac12}\\\label{eq1DnormV}
&\le C\|u-w\|_2^2\|u\|_2\|v\|_2^2,
\end{align}
and
\begin{align}\notag
&\int_{D_m}|du-dw||d^2u-d^2w||v||dv|\,dA\le\\\notag
&\|v\|_\infty\left(\int_{D_m}|d^2u-d^2w|^2\right)^{\tfrac12}
\left(\int_{D_m}|du-dv|^2|dv|^2\,dA\right)^{\tfrac12}\le\\\notag
&\|v\|_\infty\left(\int_{D_m}|d^2u-d^2w|^2\right)^{\tfrac12}
\left(\int_{D_m}|du-dv|^4\,dA\right)^{\tfrac14}
\left(\int_{D_m}|dv|^4\,dA\right)^{\tfrac14}\\\label{eq1DnormVI}
&\le C\|u-w\|_2^2\|v\|_2^2.
\end{align}

Combining \eqref{eq0Dnorm} with \eqref{eq1DnormI} --  \eqref{eq1DnormVI}, we find
\begin{equation}\label{eqdiffest}
\|F(v)-F_w(v)\|_1
\le\rho\bigl(\|u-w\|_2\bigr)\|v\|_2,
\end{equation}
where $\rho(t)$, $t\ge 0$ is a continuous function such that $\rho(t)\to 0$ as $\rho\to 0$. Taking $w=0$ in \eqref{eqdiffest}, we see that $F$ is continuous. Moreover, compactly supported smooth functions are dense in $\sblv_2(D_m,\R^n)$ and the set of compact operators is closed in the operator norm. Hence, \eqref{eqdiffest} implies that the claim follows once we prove that $F_w$ is a compact operator. To show that we estimate the $(2,2)$-norm of $F_w(v)$ using the same basic relations between Sobolev norms that was used above.

For the $L^2$-norm we have
\begin{equation}\label{eqF_w0}
\int_{D_m}|h(u)\cdot dw\cdot v|^2\,dA\le C\|v\|_2^2,
\end{equation}
since $h$ and the derivative of $w$ are bounded. For the $L^2$-norm of the derivative we have
\begin{align}\notag
&\int_{D_m}|d(h(u)\cdot dw\cdot v)|^2\,dA\le
C\int_{D_m}\bigl(|du||v|+|dv|\bigr)^2\,dA\\\label{eqF_w1}
&=C\int_{D_m}\bigl(|du|^2|v|^2 + |du||v||dv|+ |dv|^2\bigr)\,dA\le
C\bigl(1+\|u\|_2\bigr)^2\|v_2\|_2^2.
\end{align}
Finally, for the $L^2$-norm of the second derivative we have
\begin{align}\notag
\int_{D_m}&|d^2(h(u)\cdot dw\cdot v)|^2\,dA\le\\\notag
C\int_{D_m}&\bigl(|du|^2|v|+|du||dv|+|d^2v|+|d^2 u||v|\bigr)^2\,dA\le\\\notag
C\int_{D_m}
\bigl(&
|du|^4|v^2|+
|du|^2|dv|^2+
|d^2v|^2+
|d^2u|^2|v|^2\\\notag
&+|du|^3|dv||v|+
|du|^2|v||d^2 v|+
|du|^2|d^2 u||v|\\\notag
&+|du||dv||d^2 v|+
|du||dv||d^2 u||v|+
|d^2v||d^2 u||v|
\bigr)\,dA\\\notag
&\le C\Bigl(\|u\|^4_2\|v\|^2_2 +
\|u\|^3_2(\|v\|_2+\|v\|_2^2) + (\|u\|_2^2+1)\|v\|_2^2\\\label{eqF_w2} &\quad\quad+\|v\|_2^3+\|v\|_2^2\|u\|_2\Bigr),
\end{align}
where the last estimate is obtained by arguments similar to those above.
The estimates \eqref{eqF_w0} -- \eqref{eqF_w2} show that
$$
\|F_w(v)\|_2^2\le p(\|v\|_2),
$$
where $p$ is a polynomial function. In particular if $\{v_j\}$ is a sequence of functions bounded in $(2,2)$-norm then so is $\{F_w(v_j)\}$. Hence by Rellich's lemma $\{F_w(v_j)\}$ contains a subsequence
converging in the $(1,2)$-norm.
It follows that $F_w$ is a compact operator. As explained above this implies the claim and thus the lemma.
\end{proof}

Let $\Lambda_0$ and $\Lambda_1$ be (ordered) Lagrangian subspaces of $\C^n$. Define the {\em complex
angle} $\theta(\Lambda_0,\Lambda_1)\in[0,\pi)^n$ inductively as follows: If
$\dim(\Lambda_0\cap\Lambda_1)=r\ge 0$ let $\theta_1=\dots=\theta_r=0$ and let $\C^{n-r}$ denote the
Hermitian complement of $\C\otimes \Lambda_0\cap\Lambda_1$ and let
$\Lambda_i'=\Lambda_i\cap\C^{n-r}$ for $i=0,1$. If $\dim(\Lambda_0\cap\Lambda_1)=0$ then let
$\Lambda_i'=\Lambda_i$, $i=0,1$ and let $r=0$. Then $\Lambda'_0$ and $\Lambda'_1$ are Lagrangian
subspaces. Let $\alpha$ be smallest angle such that $\dim(e^{i\alpha}\Lambda_0\cap\Lambda_1)=r'>0$.
Let $\theta_{r+1}=\dots=\theta_{r+r'}=\alpha$. Now repeat the construction until $\theta_n$ has been
defined. Note that $\theta(B\Lambda_0,B\Lambda_1)=\theta(\Lambda_0,\Lambda_1)$ for every $B\in\U(n)$, where $\U(n)$ is the group of unitary $(n\times n)$-matrices,
since multiplication with $e^{i\alpha}$ commutes with everything in $\U(n)$.

For $0\le s\le 1$, let $\f_j(s)\in\U(n)$ be the matrix which in the canonical coordinates $z(j)$ is
represented by the matrix
$$
\diag(e^{-i(\pi-\theta(j)_1)s},\dots,e^{-i(\pi-\theta(j)_n)s})
$$
where $\theta(j)_i$ is the $i$-th coordinate of $\theta(j) = \theta(A^+_j, A^-_j).$ If $p$ and $q$
are consecutive punctures on $\pa D_m$ then let $I(a,b)$ denote the (oriented) path in $\pa D_m$
which connects them. Define the loop $\Gamma_A$ of Lagrangian subspaces in $\C^n$ by letting the
loop
$$
\left(A|I(p_1,p_2)\right)\ast\f_2\ast\left(A|I(p_2,p_3)\right)
\ast\f_3\ast\dots\ast\left(A|I(p_m,p_1)\right)\ast\f_1
$$
of elements of $\U(n)$ act on $\R^n\subset\C^n$.

\begin{lma}\label{fred.lma}
If $c(P,\omega)=0$ then the linearization of $\bar{\pa}_J$ at $u$ is a Fredholm operator with Fredholm index
$n+\mu(\Gamma_A).$
\end{lma}

\begin{rmk}
In the case that $c(P,\omega)\ne 0$ Lemma \ref{fred.lma} remains true except that the index equality holds only modulo $c(P,\omega)$.
\end{rmk}

\begin{proof}
Proposition 5.14 \cite{ees2} proves the exact same statement for the operator $v \mapsto
\bar{\nabla}_{J,j} v.$ The proof follows from Lemma \ref{cptK.lem} and the remark immediately
after (\ref{H1mu.eq}).
\end{proof}

%----------------------------------------------------------------
\subsection{Transversality} \label{trans.sec}
%----------------------------------------------------------------
The main result of this section is:
\begin{lma}\label{trans.lma}
Fix a chord generic Legendrian submanifold $L\subset P\times\R$ and let $(u,f,j;J)\in\cand_{2,\epsilon,\Lambda}$ be such that $\bar\pa_{J,j} u=0$.
\begin{enumerate}
\item
Let $\Lambda$ be the space of admissible almost complex around $J=J_0$ then the operator $D\bar\pa_{J_0,j}$ at $(u,f,j;J_0)$ in Lemma \ref{linearization2.lem} is surjective.
\item
Assume that $u$ is not an exceptional holomorphic disks, see Definition 6.10 \cite{ees2} and let $\Lambda$ be the space of admissible Legendrian submanifolds around $L=L_0$ then the operator $D\bar\pa_{J,j}$ at $(u,f,j;L_0)$ in Lemma \ref{linearization1.lem} is surjective.
\end{enumerate}
\end{lma}

\begin{proof}
The second statement is almost identical to Theorem 6.12 \cite{ees2} and so will not be reproved here. Consider the first statement. We ignore deformations of the conformal structure on $D_m$ since they can only help when proving
transversality. The linearization in (\ref{eqlin2}) thus reduces to
\begin{equation}
\notag
D\,\bar\pa_{J,j}[v, \lambda]=
\bar\nabla_{J,j}\,v
-\tfrac12\,J\circ (\nabla_v J)\circ\pa_{J,j}\,u
 -\tfrac12\,J\circ S_\lambda\circ\pa_{J,j}\,u
:=
\bar\nabla_{J,j}\,v + Kv +K_S.
\end{equation}

Suppose $\xi$ is an element of the cokernel of this operator. By Lemma 5.1 \cite{ees2}, we can
assume $\xi$ is $C^2$-smooth on $D_m$ and satisfies
\begin{equation}
\label{coker.eq}
\int_{D_m} \langle (\bar{\nabla}_{J,j} +K)v + K_S, \xi \rangle\, dA = 0
\end{equation}
for all $v \in T_{(u,f)} \cand_{2,\epsilon}(a {\mathbf b}).$ Here $dA$ is the area form on $D_m$ and
$\langle \cdot, \cdot \rangle$ is defined fiber wise. We wish to show that $\xi$ is zero.

Let $p$ be the positive puncture on $D_m.$ By the asymptotic behavior of holomorphic disks, see
Lemma 3.6 \cite{ees2}, we can assume that for sufficiently small $r >0$, there exists a neighborhood
$(E_0, \partial E_0) \subset (D_m, \partial D_m)$ of $p$ such that
\begin{itemize}
\item $(u(E_0), u(\partial E_0)) \subset (B(a^\ast, 2r) \Pi_P(L) \cup \partial B(a^\ast, 2r)),$
\item $\Pi_P(L) \cap \partial B(a^\ast, 2r)$ are two real analytic disjoint branches,
\item $u(\pa E_0)$ is two regular oriented curves $\gamma, \tilde{\gamma}$ on the upper and lower
branches of $\Pi_P(L) \cap \partial B(a^\ast, 2r),$ respectively.
\end{itemize}

Let $q_1, \ldots, q_r \in \pa D_m$ be the preimages under $u$
of $a^\ast$ with the property that one of the components of the punctured
neighborhood of $q_j$ in $\pa D_m$ maps to $\gamma.$
Using Taylor series as in Lemma 6.8 \cite{ees2}, we see that this set is finite.

Let $q_{r+1}, \ldots, q_s \in D_m \setminus \pa D_m$ be the preimages under $u$ of $a^\ast$ with the
property that small arc in a neighborhood of $q_j$ maps to $\gamma.$ By monotonicity, see Proposition 4.3.1
\cite{a-l}, and the maximum principle in $(B(a^\ast, 2r), J|_{B(a^\ast, 2r)}) = \C^n,$ each such
$q_j$, $r+1 \le j \le s$ accounts for at least $Cr^2$ area of the disk with $C>0$ independent of
$j.$ Thus the number of such points is finite.

For $1 \le j \le s,$ let $E_j \subset D_m$ denote the connected coordinate neighborhood
$u^{-1}(B(a^\ast, 2r))$ near $q_j.$ Let $U_1 = u(E_0)$ and $U_2$ be the Schwartz reflection of $U_1$
through $\tilde{\gamma}.$

We know from standard complex analysis that if $q \in u^{-1}(a^\ast) \setminus \{ q_1, \ldots,
q_s\}$ then $u$ must intersect $u(E_0)$ transversely at $q.$ Hence, arguing again by monotonicity
and the maximum principle, we can find $x_i \in U_i \setminus (B(a^\ast, r) \cup \Pi_P(L))$ and
small neighborhoods $B(x_i, \epsilon)$, $\epsilon \ll r,$ such that
$$
u^{-1}(B(x_i, \epsilon)) \subset \bigcup_{j=i}^s E_j.
$$
Note that $x_1 \in u(E_j)$ if and only if $x_2 \in u(E_j).$ We exclude from our list any such $j$
with $x_i \notin u(E_j).$ To simplify notation, we continue to index this possibly shortened list by
$1 \le j \le s.$

For $ 1 \le j \le r,$ we double the domain $E_j$ through its real analytic boundary $\pa E_j.$ We
also double the local map $u|_{E_j}.$ We continue to denote the open disk $E_j.$ For $0 \le j \le
s,$ let $u_j = u|_{E_j}.$ We can also double (for $1 \le j \le r$) the cokernel element $\xi$ (which is anti-holomorphic) locally and define (for $0 \le j \le s$) $\xi_j = \xi|_{E_j}.$

By Lemma 6.9 \cite{ees2}, there exists a disk $E \subset \C$ and a map $\eta$ defined on $E$ such
that for $1 \le j \le s$, there exists positive integers $k_j$ and bi-holomorphic identifications
$\phi_j$ of $E_j$ with $E$ such that $\xi_j(\phi_j(z)) = \eta(z^{k_j}).$

Via our choice of perturbations of the complex structure, we can choose $K_S$ to be supported in
$B(x_2, \epsilon).$
%Assume the support is small enough such that no point outside of $D_m \subset \bigcup_{j=i}^s E_j$
%maps into the support.
Setting $v=0$ in (\ref{coker.eq}), we get
\[
0 = \int_{D_m} \langle  K_S, \xi \rangle\, dA =
\int_{E} \langle  K_S, \sum_{j=1}^s \xi_j(\phi_j(z)) \rangle\, dA.
\]
Varying $K_S$, this implies $\sum_{j=1}^s \xi_j(\phi_j(z)) = 0.$

We can also choose $K_S$ to be supported in $B(\epsilon,x_1).$ Again setting $v=0,$ we get
\[
0 =
\int_{E} \langle  K_S, \sum_{j=1}^s \xi_j(\phi_j(z)) \rangle\, dA
+
\int_{E} \langle  K_S, \xi_0(z) \rangle\, dA
 =
\int_{E} \langle  K_S, \xi_0(z) \rangle\, dA.
\]
Varying $K_S$, this implies $\xi_0$ and hence $\xi$ vanishes.
\end{proof}

%----------------------------------------------------------------
\subsection{Gluing} \label{glue.sec}
%----------------------------------------------------------------
Let $L$ be an admissible Legendrian submanifold. Let $\M_A(a;{\mathbf b})$ and $\M_C(c;{\mathbf d})$
be moduli spaces of rigid holomorphic disks, where ${\mathbf b}$ has length $m$, $1\le j\le m$, and
${\mathbf d}$ has length $l$. Let ${\mathbf b}_{\{j\}}({\mathbf d})$ denote the word where $b_j$ is
replaced by ${\mathbf d}.$

\begin{prop}\label{glud^2=0.prop}
Assume that the $j$-th Reeb chord in ${\mathbf b}$ equals $c$. Then there exists $\delta>0$,
$\rho_0>0$ and an embedding
\begin{align}\notag
\M_A(a;{\mathbf b}) \times\M_C(c;{\mathbf d}) \times[\rho_0,\infty)&\to \M_{A+C}(a;{\mathbf
b}_{\{j\}}({\mathbf d}));\\\notag (u,v,\rho)&\mapsto u\,\sharp_\rho v,
\end{align}
such that if $u\in\M_A(a;{\mathbf b})$ and $v\in \M_C(c;{\mathbf d})$ and the image of
$w\in\M_{A+C}(a;{\mathbf b}_j({\mathbf d}))$ lies inside $B(u(D_{m+1})\cup v(D_{l+1});\delta)$ then
$w=u\,\sharp_\rho v$ for some $\rho\in[\rho_0,\infty)$.
\end{prop}

\begin{proof}
When $P \times\R = \R^{2n+1}$, this proposition appeared as Theorem 7.1 \cite{ees2}. We show that
the methods of the original proof apply to our situation.

As in most Floer gluing proofs, we show that our set-up satisfies the Floer-Picard theorem, see
Proposition 7.4 \cite{ees2} for example. Showing this has three steps: computing the
anti-holomorphicity of an "approximate glued solution"; proving the linearizations at the
approximate glued solutions are invertible in some uniform sense; and, proving a non-linear
remainder at these solutions satisfies a quadratic estimate. These three steps are presented as
Lemmas 7.5, 7.9 and 7.16 in \cite{ees2}.

For any $k$ we consider the disk $D_{k+1}$ as having a half-infinite strip $[0,1] \times [0,
\infty)$ at $p_0$ and half-infinite strips $[0,1] \times (-\infty, 0]$ at $p_j$, $1 \le j \le k.$ We
consider the disk $D_{l+m+1;\rho}$ to be gotten by removing the half-infinite strip $[0,1] \times
[\rho, \infty)$ at $p_0 \in D_{l+1}$, removing the half-infinite strip $[0,1] \times (-\infty,
-\rho]$ at $p_j \in D_{m+1}$, where $j \ne 0$, and gluing these two pieces, say $A_m$ and $A_l$
together at $\Omega = [-2,2] \times [0,1].$ This procedure also glues conformal structures.

Consider the triple $(u,v,\rho)$ from the proposition statement. Since our set-up is $\C^n$ near
$b_j^\ast,$ just as is done in Section 7.4, \cite{ees2}, we find vector fields $\xi$ and $\eta$
supported near $p_j \in D_{m+1}$ and $p_0 \in D_{l+1}$ respectively such that
\[
\exp_{b_j^\ast}(\xi(\zeta)) = u(\zeta),\quad \exp_{b_j^\ast}(\eta(\zeta)) = v(\zeta).
\]
Let $z=\tau+it$ be a complex coordinate on $\Omega$ and let $\alpha^\pm\colon\Omega\to\C$ be cut-off
functions which are real valued and holomorphic on the boundary and with $\alpha^+=1$ on
$[-2,-1]\times[0,1]$, $\alpha^+=0$ on $[0,2]\times[0,1]$, $\alpha^-=1$ on $[1,2]\times[0,1]$, and
$\alpha^-=0$ on $[-2,0]\times[0,1]$.

We then construct the approximate glued solution as follows:

\begin{equation}\notag
w_\rho(\zeta)=
\begin{cases}
&v(\zeta),\quad\quad \zeta\in A_l,\\
&u(\zeta),\quad\quad \zeta\in A_m,\\
&\exp_{b_j^\ast}(\alpha^-(z)\xi(z)+\alpha^+(z)\eta(z)),
\quad\quad z\in\Omega.
\end{cases}
\end{equation}

To show that $w_\rho$ is approximately holomorphic, we must show that
\begin{equation}\notag
\|\bar\pa w_\rho\|_1=\Ordo(e^{-\theta\rho}),
\end{equation}
where $\theta$ is the smallest component of the complex angle at the self-intersection at which the
gluing occurs. But note that away from $\Omega$, $w_\rho$ is holomorphic, while for sufficiently
large $\rho$, $w_\rho(\Omega)$ sits inside a neighborhood of $b_j^\ast$ where $(P,J) = \C^n.$ Thus
the estimate for this middle region follows from the case $P = \R^{2n}$ proved as Lemma 7.5
\cite{ees2}.

Adopting Lemma 7.9 \cite{ees2}, the uniform invertibility, to our set-up is equally straightforward.
The original proof again divides the domain into three parts as above. For sections in the domain of
$D\bar{\pa}_{J,j}|_{w\rho}$ supported at $A_l$ and $A_m$, the invertibility follows (up to
contributions from cut-off functions) from the given invertibility of $D\bar{\pa}_{J,j}|_{u}$ and
$D\bar{\pa}_{J,j}|_{v}.$ For sections supported over the middle, invertibility is verified
explicitly. As in the previous paragraph, for such sections (and sufficiently large $\rho$), our
set-up again reduces to \cite{ees2}.

It remains to prove a quadratic estimate for the non-linear term, as is done in Lemma 7.16
\cite{ees2}. Since we are proving a stationary gluing theorem, like in the previous two steps, we
fix the Legendrian submanifold $L$ and the almost complex structure $J.$ Working in a trivialization
(so that addition make sense) the non-linear term $N$ is
defined by the expansion
\begin{equation}
\label{nonlin.eq}
 \bar{\pa}_{J,j} \exp^{\sigma(\zeta)}_{w_\rho(\zeta)} \left( \tau(\zeta) \right) = \bar{\pa}_{J,j}
w_\rho (\zeta) + D\,\bar\pa_{J,j}(\tau)(\zeta) + N(\tau)(\zeta)
\end{equation}
where $\tau \in T_{w_\rho} \cand_2(a {\mathbf b}_{\{j\}} ({\mathbf d})).$ The necessary estimate to
prove is that for $\tau_1, \tau_2 \in T_{w_\rho} \cand_2(a {\mathbf b}_{\{j\}} ({\mathbf d}))$
\begin{equation}
\label{nonlinest.eq}
\|N(\tau_1)-N(\tau_2)\|_1\le C\Bigl(\|\tau_1\|_2+\|\tau_2\|_2\Bigr) \|\tau_1-\tau_2\|_2.
\end{equation}

The analogous result in the case $P=\R^{2n}$ was proved in Lemma 7.16 in \cite{ees2}. The proof has
two steps. First the function $N$ is expressed in terms of formal variables instead of actual
functions and this function is proved to have certain vanishing properties. Second the vanishing
properties and standard Sobolev estimates are combined to prove the estimate. The second step is
completely analogous in the case of $\R^{2n}\times\R$ from \cite{ees2} but the first step is somewhat
different. We will therefore describe only the first step.

In order to define the function discussed above, note that if $\nu\in T_\zeta D_m$ and if $X$ is the Jacobi field along the geodesic
\[
\gamma(s) = \exp^\sigma_{w_\rho(\zeta)} (s \tau( \zeta))
\]
satisfying the initial conditions
\begin{equation} \label{jac.eq}
X(0) = d w _\rho(\nu), \quad \nabla_{\dot{\gamma}(s)} X(0) =
\nabla_{d w_\rho (\nu)} \tau (\zeta),
\end{equation}
then $d \exp^\sigma_{w_\rho (\zeta)} \tau(\zeta) [ \nu] = X(1).$

The source of our function related to $N$ is a bundle $E\to D_{l+m+1}\times P\times\conf_{l+m+1}$, with fibers
$$
T_xP\times\mbox{Hom}(T_\zeta D_{l+m+1}, T_x P)\times\mbox{Hom}(T_xP\oplus T_xP,T_xP).
$$
We write elements in this space as $e=(x,\zeta,\kappa,y,\xi,\theta)$ where $y, \xi, \theta$ are the
fiber coordinates. To simplify notation, we ignore the variables $\sigma, \sigma'$ (on which the
metric depends) which are necessary when computing estimates, see \cite{ees2}. We will also suppress
the conformal structure variable $\kappa$ from the notation. The contributions to the non-linear
estimate from the finite-dimensional space of conformal variations are easily handled by an argument
similar to the one in \cite{ees2} Lemma 7.16 once the estimate without conformal variations has been
established. To simplify notation below we therefore consider a fixed conformal structure
$j=j_\kappa$ on $D_{m+l+1}$.

Define the bundle map $\hat{\Phi}: E \rightarrow \mbox{Hom}^{(0,1)}(T_\zeta D_{l+m+1}, T_{\exp_x(y)}P)$ by
\[
\hat{\Phi}(e) [\nu] = X[\nu] + J X[j\nu],
\]
where $X[\nu]$ is the value of the Jacobi field along the geodesic $\gamma(s)=\exp_x(sv)$, which satisfies the initial conditions
$$
X[\nu](0)=\xi[\nu],\quad \nabla_s X(0)=\theta(\xi[\nu],y),
$$
at $s=1$.

We next trivialize the bundle $\mbox{Hom}^{0,1}(T D_{l+m+1}, TP)$ in the standard way. Using notation
from the proof of Lemma~\ref{linearization1.lem}, the trivialization is a bundle map ${\Phi}: E \rightarrow
\mbox{Hom}^{(0,1)}(TD_{l+m+1}, TP)$ (which descends to the identity on $D_{l+m+1} \times P$ but which is not a fiber wise homomorphism)
\begin{equation}\label{eqdefPHI}
\Phi(e) =(1 - S_y)^{-1} \Pi_{-y} \hat{\Phi}(e).
\end{equation}

Recall
\[
\bar{\pa}_{J,j} w_\rho =d w_\rho + J d w_\rho j,
\]
and the formula for the linear term given by (\ref{eqlin1}) or (\ref{eqlin2}) (with some of the terms
disappearing because $L$ and $J$ are fixed and because we consider a fixed complex structure on $D_{m+l+1}$).
We thus define the function $\tilde{N}$ related to the non-linear term as follows
\begin{align}\notag
\tilde{N}(e)[\nu] =
&\Phi(e)[\nu] - (\xi + J \xi j)[\nu]  \\
 \label{N.eq}
&- \left(\theta(\xi[\nu],y) + J \theta(\xi[j\nu],y) - \frac{1}{2} J \nabla_y J (\xi - J \xi j)[\nu] \right),
\end{align}
where $\nu\in T_\zeta D_{m+j+1}$.

We claim that once we prove
\begin{equation}
\label{nonlinkey.eq}
\tilde{N}(x, \zeta, 0, \eta, \xi) = 0, \quad D_y \tilde{N}(x, \zeta, 0, \eta, \xi) = 0 ,
\end{equation}
the non-linear quadratic estimate (\ref{nonlinest.eq}) follows. To see this claim, note that these are exactly the equations given in
(7.95) in the proof of Lemma 7.16 \cite{ees2}. The same set of manipulations which follow (7.95) \cite{ees2} then prove a quadratic estimate for $N$.

It remains to prove (\ref{nonlinkey.eq}). To prove the first equation, we substitute $y=0$ into (\ref{eqdefPHI}). We find
$$
\Phi(e)=\xi+J\xi j.
$$
Hence,
$$
\tilde{N}(\zeta,x,0,\xi,\theta)=\xi+J\xi j-(\xi +J\xi j)=0.
$$
To show that also the derivative with respect to $y$ vanishes at $0$ we expand $\Phi(\zeta,x,\epsilon y,\xi,\theta)$ in $\epsilon$.
Arguing as in the proof of Lemma \ref{linearization1.lem} we find
\begin{equation}\label{eqagain1}
(1-S_{\epsilon y})^{-1}= 1 - \epsilon\frac12 J(\nabla_y J) +\Ordo(\epsilon^2)
\end{equation}
and
\begin{align}\label{eqagain2}
\Pi_{-\epsilon y}(\hat\Phi(e))=
\xi+J\xi j
+\epsilon\Bigl(\theta(\xi,y)+J\theta(\xi j,y)
+(\nabla_yJ)\xi j\Bigr)+\Ordo(\epsilon^2).
\end{align}
Multiplying \eqref{eqagain1} and \eqref{eqagain2} we get
\begin{align*}
\Phi(e) &= \xi + J\xi j +\epsilon\Bigl(\theta(\xi,y)+J\theta(\xi j,y)+(\nabla_yJ)\xi j-\frac12 J(\nabla_y J)(\xi + J\xi j)\Bigr) +\Ordo(\epsilon^2)\\
&=
\xi + J\xi j
+\epsilon\Bigl(\theta(\xi,y)+J\theta(\xi j,y)
- \frac{1}{2} J(\nabla_y J)(\xi-J\xi j)\Bigr)+\Ordo(\epsilon^2).
\end{align*}
Together with \eqref{N.eq} this shows that $\tilde{N}(\zeta,x,\epsilon y,\xi,\theta)=\Ordo(\epsilon^2)$ and the derivative vanishes as claimed.
\end{proof}

%----------------------------------------------------------------
\subsection{Spin structures and orientations} \label{spin.sec}
%----------------------------------------------------------------
In this Subsection, we need not make any assumption on $c(P,\omega).$
To orient the moduli spaces $\M_A(a;{\mathbf{b}})$ of $J$-holomorphic disks with boundary on a Legendrian submanifold $L\subset P\times\R$ we note (for example looking at the proof of Proposition \ref{mfld.prop}) that the tangent space of $\M_A(a;{\mathbf b})$ is closely related to the kernel of the operator $D\bar\pa_{J,j}$ in Lemmas \ref{linearization1.lem} and \ref{linearization2.lem} with $L$ respectively $J$ fixed. (In fact, if the number of Reeb chords in the word $a{\mathbf b}$ is $\ge 3$ then the tangent space of the moduli space equals this kernel and if the number is $<3$ then the tangent space is a quotient of the kernel under the natural action of the automorphism group of the source disk.) Since the linearized operators are surjective, to orient the moduli spaces it is sufficient to orient the determinant line of the operator $D\bar\pa_{J,j}$ over $\cand_{2,\epsilon}$. To accomplish this we follow the approach taken in \cite{ees3}: we orient the determinant line of the restriction $D\bar\pa_J$ of $D\bar\pa_{J,j}$ to the infinite dimensional subspace which is the complement of the space of conformal variations and orient the space of conformal variations (or automorphisms in the case that the number of punctures is smaller than $3$) separately. It is a consequence of Lemma 3.17 in \cite{ees3} that these two orientations induce an orientation on the determinant line of $D\bar\pa_{J,j}$.

For general $L$, the determinant line bundle discussed above is not necessarily orientable; however, it is if $L$ is ``relatively spin'', see \cite{FOOO}. In particular, if $L$ is spin then it is also relatively spin. We restrict attention to this case below since it is sufficient for our applications and since the construction of an orientation is completely analogous to the corresponding construction in \cite{ees3}.

Let $L\subset P\times\R$ be a Legendrian submanifold which is spin.
We do not need to make assumptions on $c(P, \omega).$
Fix an orientation and a spin structure on $L$. We view this data as giving a trivialization of the stabilized tangent bundle $\tilde TL= TL\oplus\R^2$ over the $3$-skeleton of $L$. (The fact that $L$ is spin assures that there exists a trivialization which extends over the $2$-skeleton, the extension to the $3$-skeleton is then automatic.) We pick a specific triangulation of $L$ such that all Reeb chord endpoints lie in the $0$-skeleton and such that all capping paths lie in the $1$-skeleton. We then choose a trivialization over the $3$-skeleton of this triangulation and specify capping operators at all Reeb chords as in Subsection 3.3.4 in \cite{ees3}. Also as in Subsection 3.4.2 of \cite{ees3} we consider the direct sum of the operator of interest, $D\bar\pa_J$, with the standard $\bar\pa$-operator acting on $\C^2$-valued functions on $D_m$ and with the trivialized Lagrangian boundary conditions along $\pa D_m$ defined there. We call this operator $\bar\pa_0$ and remark that it has trivial kernel and trivial cokernel. (The reason for using this stabilization (the two extra directions) is to allow for a continuous extension of the capping operators in the case of a family of Legendrian submanifolds, see the proof of \cite{ees3} Lemma 3.14.)

If $u\in\cand_{2,\epsilon}$ then we pick a complex trivialization of $(u^\ast TP)\oplus\C^2$. Homotoping $u|\pa D_m$ so that it lies in the $1$-skeleton of $L$ we get an induced trivialization of the Lagrangian boundary condition for $D\bar\pa_J\oplus\bar\pa_0$ over $u$, see Subsection 3.4.2 \cite{ees3}. Since $\tilde TL$ is trivialized also over the $2$-skeleton it is easy to see that the trivialization on the boundary of $u$ is homotopically unique. With these trivializations chosen we may view the operator $D\bar\pa_J\oplus\bar\pa_0$ as an operator
\[
\sblv_{2}[\lambda](D_m,\C^{n+2})\to
\sblv_{1}(D_m, {T^\ast}^{0,1}D_m\otimes\C^{n+2}; [0]),
\]
where the right hand side is defined in (\ref{H1mu.eq}), and using the notation of \cite{ees3},
the $[\lambda]$ on the left hand side 
indicates the subspace of the Sobolev space of functions which satisfies a Lagrangian boundary condition specified by $\lambda$ and for which the restriction of $D\bar\pa_J$ to the boundary vanishes. In \cite{ees3}, Definition 3.16, an orientation on this bundle was specified using the capping operators and the trivialization of the boundary condition discussed above. Since the choice of the complex trivialization of $TP\oplus\C^2$ is homotopically unique we get an orientation of the determinant line of $D\bar\pa_J$ over $\cand_{2,\epsilon}$ as desired. With this accomplished, the arguments that the differential in the DGA is indeed a differential and that the stable tame isomorphism class of the DGA is independent of the Legendrian isotopy class of $L$ are exactly the same as the arguments given in \cite{ees3}.

In \cite{ees3}, Theorem 4.29, the effect of changing the spin structure on the orientations is described. The situation is completely analogous here: the difference between two spin structures on $L$ is an element in $H^1(L;\Z_2)$ and the change in sign of a rigid disk is measured by the evaluation of this cohomology class on its capped off restriction to the boundary.

%----------------------------------------------------------------
\section{Applications}\label{app.sec}
%----------------------------------------------------------------
In this section we given two immediate applications of the theory described in this paper. In Subsection~\ref{knot.sec}
we describe how to construct isotopy invariants of submanifolds of $\R^n.$ Computations and applications of these
invariants will be discussed in subsequent papers. In the following subsection we discuss the question: when are two Legendrian
submanifolds in a Darboux ball which are Legendrian isotopic in the ambient manifold also Legendrian isotopic in the Darboux ball?
Specifically we show that if the Legendrian submanifolds are distinguished in the Darboux ball by contact homology then they
are also distinguished in the ambient manifold.

%----------------------------------------------------------------
\subsection{Invariants of submanifolds of $\R^n$}\label{knot.sec}
%----------------------------------------------------------------
One of the prime motivations for the work in this paper is the following application. Let $M$ be an immersed submanifold with transverse multiple points in $\R^n.$ Let $W=\R^n\times S^{n-1}$ be the unit cotangent bundle of
$\R^n.$ The restriction of the canonical 1-form $\lambda_{\R^n}$ to $W$ is a contact form on $W,$
which we denote by $\alpha.$ Now let $L_M$ be the unit conormal bundle of $M$ in $\R^n.$ That is $L_M$
consists of elements of $W$ that vanish on tangent vectors to $M.$ One may easily check that $L_M$
is a Legendrian submanifold of the contact manifold $W$ and that a regular homotopy without self tangencies of $M$ will produce a Legendrian isotopy of $L_M.$ Thus the Legendrian submanifold $L_M$ and any Legendrian isotopy invariants of it,
is invariant of $M$ up to regular homotopies without self tangencies. In particular, if $M$ is embedded then any regular homotopy of $M$ without self tangencies is and isotopy
and the Legendrian isotopy class of $L_M$ is an isotopy invariant of $M$.

To bring this situation into the setup of this paper we note that the diffeomorphism
\[
\Psi:W\to J^1(S^{n-1})= T^*S^{n-1}\times \R: ({\mathbf p}, {\mathbf q})\mapsto ({\mathbf q},
{\mathbf p} - \langle {\mathbf q}, {\mathbf p}\rangle {\mathbf q}, \langle {\mathbf q}, {\mathbf
p}\rangle)
\]
from $W$ to the 1-jet space of $S^{n-1},$ takes the contact form $\alpha$ to the contact form
$dz-\lambda_{S^{n-1}}$ described in Example~\ref{1jet.ex}. Thus we may move our Legendrian isotopy
problem from $W$ to $J^1(S^{n-1}).$ In $J^1(S^{n-1})$ the DGA defined in this paper is an invariant
of Legendrian isotopy and thus an invariant of the isotopy class of $M$ in $\R^n.$ In particular we
have:
\begin{thm}
The stable tame isomorphism class of the contact homology DGA of $L_M\subset J^1(S^{n-1})$ is an
isotopy invariant of $M$ in $\R^n.$
\end{thm}
For a leisurely introduction to the use of contact geometry in constructing invariants of
submanifolds of Euclidean space (of manifold in general) see \cite{EkholmEtnyre}.

The idea described above has been used by Ng \cite{ng1, ng2, ng3} to write down an invariant of
knots in $\R^3.$ In \cite{ng1} a knot in $\R^3$ was represented as a braid and what should be the
differential in the DGA was written down in terms of the braid. Ng proceeded to show that the
resulting DGA was in invariant of the knot using Markov's theorem that relates two braid
representations of a knot. In subsequent papers \cite{ng2, ng3} Ng has shown this DGA is a very
interesting and powerful invariant of knots.

We point out that while Ng's DGA is an interesting invariant of knots it is not obvious that it is the Contact Homology DGA. With this paper that the Contact Homology DGA is well defined for Legendrian submanifolds in a 1-jet space and a future paper will show that Ng's DGA is in fact the contact homology DGA.

%----------------------------------------------------------------
\subsection{Localizing contact homology}
%----------------------------------------------------------------
In this section we show roughly that Legendrian submanifolds in $P\times \R$ that are locally distinguished by
their contact homology are globally distinguished by their contact homology. This will allow us to
construct many non-Legendrian isotopic Legendrian submanifolds in $P\times\R.$ Recall Darboux's
theorem says that any point $p$ in a contact manifold has a neighborhood contactomorphic to a
neighborhood of the origin in $\R^{2n+1}.$ In addition, it can be shown that the contactomorphism
preserves preassigned contact forms. Thus any point in $P\times \R$ has a neighborhood $U$ that is
diffeomorphic to a neighborhood of the origin in $\R^{2n+1}$ so that $dz-\theta$ is taken to the
standard contact form on $\R^{2n+1}.$
\begin{thm}
If $L_0$ and $L_1$ are two Legendrian submanifolds of $\R^{2n+1}$ that are distinguished by their
contact homology DGA's then we may transport them to $U\subset P\times\R$ and they will still be
distinguished by their contact homology DGA's.
\end{thm}
\begin{rmk}
This theorem has bearing on the following interesting question: when are two Legendrian submanifolds in a Darboux ball which are Legendrian isotopic in the larger ambient manifold Legendrian isotopic also in the Darboux ball? To see that the answer to this question is sometimes in the negative, note that the conormal lifts of the two curves in the unit disk depicted in Figure \ref{fig:curves} are Legendrian isotopic in the unit tangent bundle of the sphere $S^2$ if the unit disk is included as a small disk around a point in $S^2$ but that they are not Legendrian isotopic in the unit tangent bundle of the disk.
\begin{figure}[htbp]
\begin{center}
\includegraphics[angle=0, width=8cm]{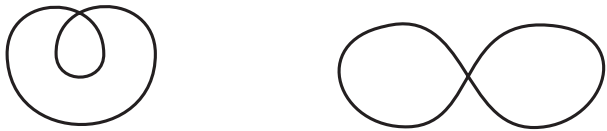}
\end{center}
\caption{Two curves.}
\label{fig:curves}
\end{figure}

\end{rmk}
\begin{proof}
By scaling $\R^{2n+1}$ we may Legendrian isotop $L_0$ and $L_1$ into a small neighborhood $V$ of the
origin that may be mapped preserving the contact form to $P\times\R.$ Since the contact forms are
preserved it is easy to see the Reeb fields are preserves. Thus the generators for the contact
homology in $\R^{2n+1}$ and in $P\times\R$ are the same. It is also clear that the projection of $V$
to $\R^{2n}$ and the projection of its image in $P\times\R$ to $P$ can be chosen to be
holomorphically equivalent (that is we can choose the almost complex structure on $P$ so that this
is the case). Moreover both the projections can be arranged so that the boundary of their closures
are $J$-convex and far from the projections of the $L_i.$ This implies, using the monotonicity lemma
for holomorphic curves, that any holomorphic curve with boundary on one of the $L_i$ must be
contained in the projection of $V.$ Thus all the holomorphic disks used in the definition of the
differential of the DGA in $\R^{2n+1}$ and in $J\times\R$ are the same. We conclude that the DGA's
are the same.
\end{proof}
Using the examples of non-Legendrian isotopic Legendrian submanifold of $\R^{2n+1}$ (with its
standard contact structure) constructed in \cite{ees1, ees3} we have the following immediate
corollary.
\begin{cor}
In any $P\times\R$ there are infinitely many Legendrian submanifold that are not Legendrian isotopic
but have the same classical invariants (that is topological isotopy type, Thurston-Bennequin
invariant and rotation class).
\end{cor}

\end{document}